\documentclass[10pt]{article}
\usepackage{setspace}

\usepackage{amsmath} 
\usepackage{amssymb}
\usepackage{latexsym}
\usepackage{xcolor}
\usepackage{fancyhdr}
\allowdisplaybreaks 

 \usepackage{comment}

\def\XXint#1#2#3{{\setbox0=\hbox{$#1{#2#3}{\int}$} 
\vcenter{\hbox{$#2#3$}}\kern-.5\wd0}}   

 \numberwithin{equation}{section}
\newtheorem{theorem}[equation]{Theorem}
\newtheorem{proposition}[equation]{Proposition}
\newtheorem{definition}[equation]{Definition}
\newtheorem{remark}[equation]{Remark}
\newtheorem{lemma}[equation]{Lemma}
\newtheorem{corollary}[equation]{Corollary}

\title{
A nonvariational Neumann problem  for the Helmholtz equation
 } 
 
\author{  
Massimo Lanza de Cristoforis
\\
Dipartimento di Matematica `Tullio Levi-Civita', 
\\
Universit\`a degli Studi di Padova, 
\\
Via Trieste 63, Padova 35121, 
Italy. 
\\
E-mail: mldc@math.unipd.it   }

\date{\ }

\begin{document}
 
 \maketitle

\noindent
{\bf Abstract:}  We consider a bounded open subset $\Omega$ of ${\mathbb{R}}^n$ of class $C^{1,\alpha}$ for some $\alpha\in]0,1[$ and we solve the Neumann problem for the Helmholtz equation both in $\Omega$ and  in the exterior  of $\Omega$. We look for  solutions in the space for $\alpha$-H\"{o}lder continuous functions  that may   not have a classical normal derivative at the boundary points of $\Omega$ and that may have an infinite Dirichlet integral around the boundary of $\Omega$. Namely for solutions  that  do not belong to the classical variational setting. 

 \vspace{\baselineskip}

\footnotetext[1]{
{{\bf 2020 Mathematics Subject Classification:}}   31B10, 
35J25, 35J05.
}
\footnotetext[2]{
{\bf Keywords:}   H\"{o}lder continuity,  Helmholtz equation, Schauder spaces with negative exponent, nonvariational solutions. 
}
\section{Introduction} Our starting point are the classical examples of Prym \cite{Pr1871} and  Hadamard \cite{Ha1906} of harmonic functions in the unit ball of the plane that are continuous up to the boundary and have  infinite Dirichlet integral, \textit{i.e.},  whose
 gradient is not square-summable. Such functions solve the classical Dirichlet problem in the unit ball, but not the corresponding weak (variational) formulation. For a discussion on this point we refer to 
Maz’ya and Shaposnikova \cite{MaSh98}, Bottazzini and Gray \cite{BoGr13} and
 Bramati,  Dalla Riva and  Luczak~\cite{BrDaLu23} which contains examples of H\"{o}lder continuous harmonic functions with infinite Dirichlet integral.

In the papers \cite{La24b}, \cite{La24c}, the author has analyzed the Neumann problem for the Laplace and the Poisson equation. In particular, the author has introduced a distributional normal derivative on the boundary  for $\alpha$-H\"{o}lder continuous functions that have a Laplace operator in a Schauder space of negative exponent and has characterized the corresponding space for the first order traces.

In  \cite{La25} the author has considered a corresponding outward distributional normal derivative on the boundary and has  the analyzed the  uniqueness problem for  boundary value problems for the Helmholtz equation in exterior domains with $\alpha$-H\"{o}lder continuous solutions    that may have infinite Dirichlet integral around the boundary points  and that satisfy the Sommerfeld radiation condition. In \cite{La25a} the author has proved those
  potential theoretic properties of the single layer potential with distributional densities in the above mentioned space of first order traces of  $\alpha$-H\"{o}lder continuous functions that are necessary for the potential theoretic treatment of boundary value problems for the Helmholtz equation with nonvariational solutions
  
  In the present paper we prove the classical result on the validity of the Sommerfeld radiation condition for single layer potentials with distributional density in the above mentioned space of first order traces of  $\alpha$-H\"{o}lder continuous functions (cf.~Theorem \ref{thm:pora}),  we prove a third Green Identity for  $\alpha$-H\"{o}lder continuous functions that have a Laplace operator in a Schauder space of negative exponent for interior and exterior domains (cf.~Theorems \ref{thm:dthirdgreenh}, \ref{thm:dethirdgreenh}), we solve the interior    and the exterior  Neumann problem with distributional data in the above mentioned space of first order traces of  $\alpha$-H\"{o}lder continuous functions (cf.~Theorems \ref{thm:exintneuhe}, \ref{thm:exextneuhe}). Especially for the exterior problem, we note that we also consider the case in which the domain is multiply connected.  
  
  The paper is organized as follows. Section \ref{sec:prelnot} is a section of
 preliminaries and notation. In section \ref{sec:doun}, we introduce the distributional form of the outward unit normal derivative of  \cite{La25} for functions which are locally  $\alpha$-H\"{o}lder continuous and that have  locally a Laplace operator in a Schauder space of negative exponent. In section \ref{sec:prefusoh}, we introduce some preliminaries on the fundamental solution of the Helmholtz equation. In section \ref{sec:acsilah-1a}, we introduce some preliminaries on the acoustic layer potentials with densities in the above mentioned space of first order traces of  $\alpha$-H\"{o}lder continuous functions. In section \ref{sec:laporaco}, we discuss the validity of the Sommerfeld radiation condition for 
   (distributional) acoustic  layer potentials. In section \ref{sec:acvolpot}, we  introduce some preliminaries on the acoustic volume potential with density in a Schauder space with either positive or negative exponent. In section \ref{sec:dgridh}, we prove  third Green Identity for the Helmholtz equation for $\alpha$-H\"{o}lder continuous functions that have a Laplace operator in a Schauder space of negative exponent. In section \ref{sec:rinefh}, we introduce a known   representation theorem for the interior Neumann eigenfunctions of the Laplace operator. In section \ref{sec:inpddh}, we solve the interior Neumann problem for the Helmholtz equation with distributional data in the above mentioned space of first order traces of  $\alpha$-H\"{o}lder continuous functions.  In section \ref{sec:renefh}, we introduce a known   representation theorem for the exterior Neumann eigenfunctions of the Laplace operator. In section \ref{sec:enpddh}, we solve the exterior Neumann problem for the Helmholtz equation with distributional data in the above mentioned space of first order traces of  $\alpha$-H\"{o}lder continuous functions.

 \section{Preliminaries and notation}\label{sec:prelnot} Unless otherwise specified,  we assume  throughout the paper that
\[
n\in {\mathbb{N}}\setminus\{0,1\}\,,
\]
where ${\mathbb{N}}$ denotes the set of natural numbers including $0$. 
If $X$ and $Y$, $Z$ are normed spaces, then ${\mathcal{L}}(X,Y)$ denotes the space of linear and continuous maps from $X$ to $Y$ and ${\mathcal{L}}^{(2)}(X\times Y, Z)$ 
denotes the space of bilinear and continuous maps from $X\times Y$ to $Z$ with their usual operator norm (cf.~\textit{e.g.}, \cite[pp.~16, 621]{DaLaMu21}). 
 
$|A|$ denotes the operator norm of a matrix $A$ with real (or complex) entries, 
       $A^{t}$ denotes the transpose matrix of $A$.	 
 
 Let $\Omega$ be an open subset of ${\mathbb{R}}^n$. $C^{1}(\Omega)$ denotes the set of continuously differentiable functions from $\Omega$ to ${\mathbb{C}}$. 
 Let $s\in {\mathbb{N}}\setminus\{0\}$, $f\in \left(C^{1}(\Omega)\right)^{s} $. Then   $Df$ denotes the Jacobian matrix of $f$.   
 
 For the (classical) definition of   open Lipschitz subset of ${\mathbb{R}}^n$ and of   open subset of ${\mathbb{R}}^n$ 
   of class $C^{m}$ or of class $C^{m,\alpha}$
  and of the H\"{o}lder and Schauder spaces $C^{m,\alpha}(\overline{\Omega})$
  on the closure $\overline{\Omega}$ of  an open set $\Omega$ and 
  of the H\"{o}lder and Schauder spaces
   $C^{m,\alpha}(\partial\Omega)$ 
on the boundary $\partial\Omega$ of an open set $\Omega$ for some $m\in{\mathbb{N}}$, $\alpha\in ]0,1]$, we refer for example to
    Dalla Riva, the author and Musolino  \cite[\S 2.3, \S 2.6, \S 2.7, \S 2.9, \S 2.11, \S 2.13,   \S 2.20]{DaLaMu21}.  If $m\in {\mathbb{N}}$, 
 $C^{m}_b(\overline{\Omega})$ denotes the space of $m$-times continuously differentiable functions from $\Omega$ to ${\mathbb{C}}$ such that all 
the partial derivatives up to order $m$ have a bounded continuous extension to    $\overline{\Omega}$ and we set
\[
\|f\|_{   C^{m}_{b}(
\overline{\Omega} )   }\equiv
\sum_{|\eta|\leq m}\, \sup_{x\in \overline{\Omega}}|D^{\eta}f(x)|
\qquad\forall f\in C^{m}_{b}(
\overline{\Omega} )\,.
\]
If $\alpha\in ]0,1]$, then 
$C^{m,\alpha}_b(\overline{\Omega})$ denotes the space of functions of $C^{m}_{b}(
\overline{\Omega}) $  such that the  partial derivatives of order $m$ are $\alpha$-H\"{o}lder continuous in $\Omega$. Then we equip $C^{m,\alpha}_{b}(\overline{\Omega})$ with the norm
\[
\|f\|_{  C^{m,\alpha}_{b}(\overline{\Omega})  }\equiv 
\|f\|_{  C^{m }_{b}(\overline{\Omega})  }
+\sum_{|\eta|=m}|D^{\eta}f|_{\alpha}\qquad\forall f\in C^{m,\alpha}_{b}(\overline{\Omega})\,,
\]
where $|D^{\eta}f|_{\alpha}$ denotes the $\alpha$-H\"{o}lder constant of the partial derivative $D^{\eta}f$ of order $\eta$ of $f$ in $\Omega$. If $\Omega$ is bounded, we obviously have $C^{m }_{b}(\overline{\Omega})=C^{m } (\overline{\Omega})$ and $C^{m,\alpha}_{b}(\overline{\Omega})=C^{m,\alpha} (\overline{\Omega})$. 
  Then $C^{m,\alpha}_{{\mathrm{loc}}}(\overline{\Omega }) $  denotes 
the space  of those functions $f\in C^{m}(\overline{\Omega} ) $ such that $f_{|\overline{\Omega'}} $ belongs to $
C^{m,\alpha}(   \overline{ \Omega' }   )$ for all bounded open subsets $\Omega'$ of ${\mathbb{R}}^n$ such that $\overline{\Omega'}\subseteq\overline{\Omega}$. 
The space of complex valued functions of class $C^m$ with compact support in an open set $\Omega$ of ${\mathbb{R}}^n$ is denoted $C^m_c(\Omega)$ and similarly for $C^\infty_c(\Omega)$. We also set ${\mathcal{D}}(\Omega)\equiv C^\infty_c(\Omega)$. Then the dual ${\mathcal{D}}'(\Omega)$ is known to be the space of distributions in $\Omega$. The support either of a function or of a distribution is denoted by the abbreviation `${\mathrm{supp}}$'.  We also set
\begin{equation}\label{eq:exto}
\Omega^-\equiv {\mathbb{R}}^n\setminus\overline{\Omega}\,.
\end{equation} 
 If $\Omega$ is a bounded open Lipschitz subset of ${\mathbb{R}}^n$, then $\Omega$  is known to have has at most a finite number of pairwise disjoint  connected components, which we denote by $\Omega_{1}$,\dots, $\Omega_{\varkappa^{+}}$ and which are open and $\Omega^-$ is known to have  at most a finite number of  pairwise disjoint  connected components that are open, which we denote by $(\Omega^{-})_{0}$,
$(\Omega^{-})_{1}$, \dots, $(\Omega^{-})_{\varkappa^{-}}$.   One and only one of such connected components is unbounded. We denote it by $(\Omega^{-})_{0}$  (cf.~\textit{e.g.}, \cite[Lemma~2.38]{DaLaMu21}).\par  

We denote by $\nu_\Omega$ or simply by $\nu$ the outward unit normal of $\Omega$ on $\partial\Omega$. Then $\nu_{\Omega^-}=-\nu_\Omega$ is the outward unit normal of $\Omega^-$ on $\partial\Omega=\partial\Omega^-$.

We now summarize the definition and some elementary properties of the Schau\-der space $C^{-1,\alpha}(\overline{\Omega})$ 
by following the presentation of Dalla Riva, the author and Musolino \cite[\S 2.22]{DaLaMu21}.
\begin{definition} 
\label{defn:sch-1}\index{Schauder space!with negative exponent}
 Let $\alpha\in]0,1]$. Let $\Omega$ be a bounded open subset of ${\mathbb{R}}^{n}$. We denote by $C^{-1,\alpha}(\overline{\Omega})$ the subspace 
 \[
 \left\{
 f_{0}+\sum_{j=1}^{n}\frac{\partial}{\partial x_{j}}f_{j}:\,f_{j}\in 
 C^{0,\alpha}(\overline{\Omega})\ \forall j\in\{0,\dots,n\}
 \right\}\,,
 \]
 of the space of distributions ${\mathcal{D}}'(\Omega)$  in $\Omega$ and we set
 \begin{eqnarray}
\label{defn:sch-2}
\lefteqn{
\|f\|_{  C^{-1,\alpha}(\overline{\Omega})  }
\equiv\inf\biggl\{\biggr.
\sum_{j=0}^{n}\|f_{j}\|_{ C^{0,\alpha}(\overline{\Omega})  }
:\,
}
\\ \nonumber
&&\qquad\qquad\qquad\qquad
f=f_{0}+\sum_{j=1}^{n}\frac{\partial}{\partial x_{j}}f_{j}\,,\ 
f_{j}\in C^{0,\alpha}(\overline{\Omega})\ \forall j\in \{0,\dots,n\}
\biggl.\biggr\}\,.
\end{eqnarray}
\end{definition}
$(C^{-1,\alpha}(\overline{\Omega}), \|\cdot\|_{  C^{-1,\alpha}(\overline{\Omega})  })$ is known to be a Banach space and  is continuously embedded into ${\mathcal{D}}'(\Omega)$. Also, the definition of the norm $\|\cdot\|_{  C^{-1,\alpha}(\overline{\Omega})  }$ implies that $C^{0,\alpha}(\overline{\Omega})$ is continuously embedded into $C^{-1,\alpha}(\overline{\Omega})$ and that the partial differentiation $\frac{\partial}{\partial x_{j}}$ is continuous from 
$C^{0,\alpha}(\overline{\Omega})$ to $C^{-1,\alpha}(\overline{\Omega})$ for all $j\in\{1,\dots,n\}$.  Generically, the  elements of $C^{-1,\alpha}(\overline{\Omega})$ for $\alpha\in]0,1[$ are not integrable functions, but distributions in $\Omega$.  Then we have the following statement of \cite[Prop.~3.1]{La24c} that shows that the elements of $C^{-1,\alpha}(\overline{\Omega}) $, which belong to the dual of ${\mathcal{D}}(\Omega)$, can be extended to elements of the dual of $C^{1,\alpha}(\overline{\Omega})$.  
\begin{proposition}\label{prop:nschext}
 Let $\alpha\in]0,1[$. Let $\Omega$ be a bounded open Lipschitz subset of ${\mathbb{R}}^{n}$.  There exists one and only one  linear and continuous extension operator $E^\sharp_\Omega$ from $C^{-1,\alpha}(\overline{\Omega})$ to $\left(C^{1,\alpha}(\overline{\Omega})\right)'$ such that
 \begin{eqnarray}\label{prop:nschext2}
\lefteqn{
\langle E^\sharp_\Omega[f],v\rangle 
}
\\ \nonumber
&&\ \
 =
\int_{\Omega}f_{0}v\,dx+\int_{\partial\Omega}\sum_{j=1}^{n} (\nu_{\Omega})_{j}f_{j}v\,d\sigma
 -\sum_{j=1}^{n}\int_{\Omega}f_{j}\frac{\partial v}{\partial x_j}\,dx
\quad \forall v\in C^{1,\alpha}(\overline{\Omega})
\end{eqnarray}
for all $f=  f_{0}+\sum_{j=1}^{n}\frac{\partial}{\partial x_{j}}f_{j}\in C^{-1,\alpha}(\overline{\Omega}) $. Moreover, 
\begin{equation}\label{prop:nschext1}
E^\sharp_\Omega[f]_{|\Omega}=f\,, \ i.e.,\ 
\langle E^\sharp_\Omega[f],v\rangle =\langle f,v\rangle \qquad\forall v\in {\mathcal{D}}(\Omega)
\end{equation}
for all $f\in C^{-1,\alpha}(\overline{\Omega})$ and
\begin{equation}\label{prop:nschext3}
\langle E^\sharp_\Omega[f],v\rangle =\langle f,v\rangle \qquad\forall v\in C^{1,\alpha}(\overline{\Omega})
\end{equation}
for all $f\in C^{0,\alpha}(\overline{\Omega})$.
\end{proposition}
When no ambiguity can arise, we simply write $E^\sharp$ instead of $E^\sharp_\Omega$. To see why the extension operator $E^\sharp$ can be considered as `canonical', we refer to \cite[\S 2]{La24d}. Next we introduce the following multiplication lemma. For a proof, we refer to \cite[\S 2]{La25}
\begin{lemma}\label{lem:multc1ac-1a}
  Let   $\alpha\in ]0,1[$. Let $\Omega$ be a bounded open Lipschitz subset of ${\mathbb{R}}^{n}$. Then the pointwise product is bilinear and continuous from 
  $C^{1,\alpha}(\overline{\Omega})\times C^{-1,\alpha}(\overline{\Omega})$ to $C^{-1,\alpha}(\overline{\Omega})$. 
\end{lemma}
 Next we introduce our space for the solutions.

\begin{definition}\label{defn:c0ade}
 Let   $\alpha\in ]0,1]$. Let $\Omega$ be a bounded open  subset of ${\mathbb{R}}^{n}$. Let
 \begin{eqnarray}\label{defn:c0ade1}
C^{0,\alpha}(\overline{\Omega})_\Delta
&\equiv&\biggl\{u\in C^{0,\alpha}(\overline{\Omega}):\,\Delta u\in C^{-1,\alpha}(\overline{\Omega})\biggr\}\,,
\\ \nonumber
\|u\|_{ C^{0,\alpha}(\overline{\Omega})_\Delta }
&\equiv& \|u\|_{ C^{0,\alpha}(\overline{\Omega})}
+\|\Delta u\|_{C^{-1,\alpha}(\overline{\Omega})}
\qquad\forall u\in C^{0,\alpha}(\overline{\Omega})_\Delta\,.
\end{eqnarray}
\end{definition}
Since $C^{0,\alpha}(\overline{\Omega})$ and $C^{-1,\alpha}(\overline{\Omega})$ are Banach spaces,   $\left(\|u\|_{ C^{0,\alpha}(\overline{\Omega})_\Delta }, \|\cdot \|_{ C^{0,\alpha}(\overline{\Omega})_\Delta }\right)$ is a Banach space.  For subsets $\Omega$ that are not necessarily bounded, we introduce the following statement of 
\cite[\S 2]{La25}.

\begin{lemma}\label{lem:c1alcof}
 Let   $\alpha\in ]0,1]	$. Let $\Omega$ be an open  subset of ${\mathbb{R}}^{n}$. Then the space
\begin{eqnarray}\label{lem:c1alcof1}
\lefteqn{
 C^{0,\alpha}_{	{\mathrm{loc}}	}(\overline{\Omega})_\Delta\equiv\biggl\{
 f\in C^{0}(\overline{\Omega}):\, f_{|\overline{\Omega'}}\in C^{0,\alpha}(\overline{\Omega'})_\Delta\ \text{for\ all\  } 
 }
\\ \nonumber
&&\qquad\qquad\qquad\qquad 
 \text{bounded\ open\ subsets}\ \Omega'\ \text{of}\ {\mathbb{R}}^n\ \text{such\ that} \  \overline{\Omega'}\subseteq\overline{\Omega}  \biggr\}\,,
\end{eqnarray}
 with the family of seminorms
\begin{eqnarray}\label{lem:c1alcof2}
\lefteqn{
{\mathcal{P}}_{C^{0,\alpha}_{	{\mathrm{loc}}	}(\overline{\Omega})_\Delta}\equiv
\biggl\{
\|\cdot\|_{C^{0,\alpha}(\overline{\Omega'})_\Delta}:\,  
}
\\ \nonumber
&& \qquad\qquad\qquad\qquad 
 \Omega'\ \text{is\ a\ bounded\ open  subset\ of\ }  {\mathbb{R}}^n \ \text{such\ that} \  \overline{\Omega'}\subseteq\overline{\Omega}
\biggr\}
\end{eqnarray}
is a Fr\'{e}chet space.
\end{lemma}
Next we introduce the following multiplication lemma. For a proof, we refer to \cite[\S 2]{La25}
\begin{lemma}\label{lem:multc2ac0ad}
  Let   $\alpha\in ]0,1[$. Let $\Omega$ be a bounded open Lipschitz subset of ${\mathbb{R}}^{n}$. Then the pointwise product is bilinear and continuous from 
  $C^{2,\alpha}(\overline{\Omega})\times C^{0,\alpha}(\overline{\Omega})_\Delta$ to $C^{0,\alpha}(\overline{\Omega})_\Delta$. 
\end{lemma}
 We also note that the following elementary lemma holds. For a proof, we refer to \cite[\S 2]{La25}
\begin{lemma}\label{lem:c1alf}
 Let   $\alpha\in ]0,1[$. Let $\Omega$ be an open  subset of ${\mathbb{R}}^{n}$. Then the space
\begin{eqnarray}\label{lem:c1alf1}
\lefteqn{
 C^{1,\alpha}_{	{\mathrm{loc}}	}(\Omega)\equiv\biggl\{
 f\in C^{1}(\Omega):\, f_{|\overline{\Omega'}}\in C^{1,\alpha}(\overline{\Omega'})\ \text{for\ all}\  \Omega'\ \text{such\ that}\ 
 }
\\ \nonumber
&&\qquad\qquad 
 \Omega'\ \text{is\ a\ bounded\ open  subset \ of}\ {\mathbb{R}}^n \,, \overline{\Omega'}\subseteq\Omega  \biggr\}\,,
\end{eqnarray}
 with the family of seminorms
\begin{eqnarray}\label{lem:c1alf2}
\lefteqn{
{\mathcal{P}}_{C^{1,\alpha}_{{\mathrm{loc}}}(\Omega)}\equiv
\biggl\{
\|\cdot\|_{C^{1,\alpha}(\overline{\Omega'})}:\,  
}
\\ \nonumber
&& \qquad\qquad 
 \Omega'\ \text{is\ a\ bounded\ open  \ subset \ of} \ {\mathbb{R}}^n \,, \overline{\Omega'}\subseteq\Omega 
\biggr\}
\end{eqnarray}
is a Fr\'{e}chet space.
\end{lemma}
 Then we can prove the following. For a proof, we refer to \cite[\S 2]{La25}.
\begin{lemma}\label{lem:c0ademc1a}
  Let   $\alpha\in ]0,1[$. Let $\Omega$ be a bounded open  subset of ${\mathbb{R}}^{n}$. Then $C^{0,\alpha}(\overline{\Omega})_\Delta$ is continuously embedded into the 
  Fr\'{e}chet space $C^{1,\alpha}_{{\mathrm{loc}}}(\Omega)$ with the family of seminorms ${\mathcal{P}}_{C^{1,\alpha}_{{\mathrm{loc}}}(\Omega)}$.
\end{lemma}
 Next we   introduce a subspace of  $C^{1,\alpha}(\overline{\Omega})$ that we need in the sequel.  
\begin{definition}\label{defn:c1ade}
 Let   $\alpha\in ]0,1[$. Let $\Omega$ be a bounded open  subset of ${\mathbb{R}}^{n}$. Let
 \begin{eqnarray}\label{defn:c1ade1}
C^{1,\alpha}(\overline{\Omega})_\Delta
&\equiv&\biggl\{u\in C^{1,\alpha}(\overline{\Omega}):\,\Delta u\in C^{0,\alpha}(\overline{\Omega})\biggr\}\,,
\\ \nonumber
\|u\|_{ C^{1,\alpha}(\overline{\Omega})_\Delta }
&\equiv& \|u\|_{ C^{1,\alpha}(\overline{\Omega})}
+\|\Delta u\|_{C^{0,\alpha}(\overline{\Omega})}
\qquad\forall u\in C^{1,\alpha}(\overline{\Omega})_\Delta\,.
\end{eqnarray}
\end{definition}
If $\Omega$ is a bounded open subset of ${\mathbb{R}}^{n}$, then 
 $C^{1,\alpha}(\overline{\Omega})$ and $C^{0,\alpha}(\overline{\Omega})$ are Banach spaces and accordingly   $\left(\|u\|_{ C^{1,\alpha}(\overline{\Omega})_\Delta }, \|\cdot \|_{ C^{1,\alpha}(\overline{\Omega})_\Delta }\right)$ is a Banach space. We also note that  if $\Omega$ is a bounded open Lipschitz subset of ${\mathbb{R}}^n$, then 
$C^{1,\alpha}(\overline{\Omega})$ is continuously embedded into $ C^{0,\alpha}(\overline{\Omega})$ and accordingly $C^{1,\alpha}(\overline{\Omega})_\Delta$ is continuously embedded into $ C^{0,\alpha}(\overline{\Omega})_\Delta$. Moreover, by applying \cite[Lem.~5.12]{La24c} to each closed ball  contained in $\Omega$, we deduce that
\begin{equation}\label{eq:c1adc2}
C^{1,\alpha}(\overline{\Omega})_\Delta\subseteq C^2(\Omega)\,.
\end{equation}
We also note that  if $\Omega$ is a bounded open Lipschitz subset of ${\mathbb{R}}^n$, then 
$C^{1,\alpha}(\overline{\Omega})_\Delta$ is continuously embedded into $ C^{0,\alpha}(\overline{\Omega})_\Delta$.  
 Next we introduce the following classical result on the Green operator for the interior Dirichlet problem. For a proof, we refer for example to \cite[Thm.~4.8]{La24b}.
 \begin{theorem}\label{thm:idwp}
Let $m\in {\mathbb{N}}$, $\alpha\in ]0,1[$. Let $\Omega$ be a bounded open  subset of ${\mathbb{R}}^{n}$ of class $C^{\max\{m,1\},\alpha}$.
Then the map ${\mathcal{G}}_{\Omega,d,+}$  from $C^{m,\alpha}(\partial\Omega)$ to the closed subspace 
 \begin{equation}\label{thm:idwp1}
C^{m,\alpha}_h(\overline{\Omega}) \equiv \{
u\in C^{m,\alpha}(\overline{\Omega}), u\ \text{is\ harmonic\ in}\ \Omega\}
\end{equation}
of $ C^{m,\alpha}(\overline{\Omega})$ that takes $v$ to the only solution $v^\sharp$ of the Dirichlet problem
\begin{equation}\label{defn:cinspo3}
\left\{
\begin{array}{ll}
 \Delta v^\sharp=0 & \text{in}\ \Omega\,,
 \\
v^\sharp_{|\partial\Omega} =v& \text{on}\ \partial\Omega 
\end{array}
\right.
\end{equation}
is a linear homeomorphism.
\end{theorem}
Next we introduce the following two approximation lemmas.
\begin{lemma}\label{lem:apr0ad}
 Let $\alpha\in]0,1[$. 
 Let $\Omega$ be a bounded open subset of ${\mathbb{R}}^n$ of class $C^{1,\alpha}$.  If $u\in C^{0,\alpha}(\overline{\Omega})_\Delta$, then there exists a sequence $\{u_j\}_{j\in {\mathbb{N}}}$ in 
 $C^{1,\alpha}(\overline{\Omega})_\Delta$ such that
 \begin{equation}\label{lem:apr1ad1}
 \sup_{j\in {\mathbb{N}}}\|u_j\|_{
 C^{0,\alpha}(\overline{\Omega})_\Delta
 }<+\infty\,,\qquad
 \lim_{j\to\infty}u_j=u\quad\text{in}\ C^{0,\beta}(\overline{\Omega})_\Delta\quad \forall\beta\in]0,\alpha[\,.
 \end{equation}
\end{lemma}
For a proof, we refer to \cite[Lem.~5.14]{La24c}.
\begin{lemma}\label{lem:apr1a}
 Let $\alpha\in]0,1]$. 
 Let $\Omega$ be a bounded open subset of ${\mathbb{R}}^n$ of class $C^{1,\alpha}$.  If $g\in C^{1,\alpha}(\overline{\Omega})$, then there exists a sequence $\{g_j\}_{j\in {\mathbb{N}}}$ in 
 $C^{\infty}(\overline{\Omega})$ such that
 \begin{equation}\label{lem:apr1a1}
 \sup_{j\in {\mathbb{N}}}\|g_j\|_{
 C^{1,\alpha}(\overline{\Omega}) 
 }<+\infty\,,\qquad
 \lim_{j\to\infty}g_j=g\quad\text{in}\ C^{1,\beta}(\overline{\Omega})\quad \forall\beta\in]0,\alpha[\,.
 \end{equation}
\end{lemma}
For a proof, we refer to \cite[Lem.~A.3]{La24c}.   Next we plan to introduce the normal derivative of the functions in  $C^{0,\alpha}(\overline{\Omega})_\Delta$ as in \cite{La24c}. To do so,    we introduce the (classical) interior Steklov-Poincar\'{e} operator (or interior  Dirichlet-to-Neumann map).
\begin{definition}\label{defn:cinspo}
 Let $\alpha\in]0,1[$.  Let  $\Omega$ be a  bounded open subset of ${\mathbb{R}}^{n}$ of class $C^{1,\alpha}$. The classical interior Steklov-Poincar\'{e} operator is defined to be the operator $S_{\Omega,+}$ from
 \begin{equation}\label{defn:cinspo1}
C^{1,\alpha}(\partial\Omega)\quad\text{to}\quad C^{0,\alpha}(\partial\Omega)
\end{equation}
that takes $v\in C^{1,\alpha}(\partial\Omega)$ to the function 
 \begin{equation}\label{defn:cinspo2}
S_{\Omega,+}[v](x)\equiv \frac{\partial  }{\partial\nu}{\mathcal{G}}_{\Omega,d,+}[v](x)\qquad\forall x\in\partial\Omega\,.
\end{equation}
 \end{definition}
   Since   the classical normal derivative is continuous from $C^{1,\alpha}(\overline{\Omega})$ to $C^{0,\alpha}(\partial\Omega)$, the continuity of ${\mathcal{G}}_{\Omega,d,+}$ implies  that $S_{\Omega,+}[\cdot]$ is linear and continuous from 
  $C^{1,\alpha}(\partial\Omega)$ to $C^{0,\alpha}(\partial\Omega)$. Then we have the following definition of \cite[(41)]{La24c}.
  \begin{definition}\label{defn:conoderdedu}
 Let $\alpha\in]0,1[$.  Let  $\Omega$ be a  bounded open subset of ${\mathbb{R}}^{n}$ of class $C^{1,\alpha}$. If  $u\in C^{0}(\overline{\Omega})$ and $\Delta u\in  C^{-1,\alpha}(\overline{\Omega})$, then we define the distributional  normal derivative $\partial_\nu u$
 of $u$ to be the only element of the dual $(C^{1,\alpha}(\partial\Omega))'$ that satisfies the following equality
 \begin{equation}\label{defn:conoderdedu1}
\langle \partial_\nu u ,v\rangle \equiv\int_{\partial\Omega}uS_{\Omega,+}[v]\,d\sigma
+\langle E^\sharp_\Omega[\Delta u],{\mathcal{G}}_{\Omega,d,+}[v]\rangle 
\qquad\forall v\in C^{1,\alpha}(\partial\Omega)\,.
\end{equation}
\end{definition}
 The normal derivative of Definition \ref{defn:conoderdedu} extends the classical one in the sense that if $u\in C^{1,\alpha}(\overline{\Omega})$, then under the assumptions on $\alpha$ and $\Omega$  of Definition \ref{defn:conoderdedu}, we have 
 \begin{equation}\label{lem:conoderdeducl1}
 \langle \partial_\nu u ,v\rangle \equiv\int_{\partial\Omega}\frac{\partial u}{\partial\nu}v\,d\sigma
  \quad\forall v\in C^{1,\alpha}(\partial\Omega)\,,
\end{equation}
where $\frac{\partial u}{\partial\nu}$ in the right hand side denotes the classical normal derivative of $u$ on $\partial\Omega$ (cf.~\cite[Lem.~5.5]{La24c}). In the sequel, we use the classical symbol $\frac{\partial u}{\partial\nu}$ also for  $\partial_\nu u$ when no ambiguity can arise.\par  
  
Next we introduce the  function space $V^{-1,\alpha}(\partial\Omega)$ on the boundary of $\Omega$ for the normal derivatives of the functions of $C^{0,\alpha}(\overline{\Omega})_\Delta$ as in  \cite[Defn.~13.2, 15.10, Thm.~18.1]{La24b}.
\begin{definition}\label{defn:v-1a}
Let   $\alpha\in ]0,1[$. Let $\Omega$ be a bounded open  subset of ${\mathbb{R}}^{n}$ of class $C^{1,\alpha}$. Let 
\begin{eqnarray}\label{defn:v-1a1}
 \lefteqn{V^{-1,\alpha}(\partial\Omega)\equiv \biggl\{\mu_0+S_{\Omega,+}^t[\mu_1]:\,\mu_0, \mu_1\in C^{0,\alpha}(\partial\Omega)
\biggr\}\,,
}
\\ \nonumber
\lefteqn{
\|\tau\|_{  V^{-1,\alpha}(\partial\Omega) }
\equiv\inf\biggl\{\biggr.
 \|\mu_0\|_{ C^{0,\alpha}(\partial\Omega)  }+\|\mu_1\|_{ C^{0,\alpha}(\partial\Omega)  }
:\,
 \tau=\mu_0+S_{\Omega,+}^t[\mu_1]\biggl.\biggr\}\,,
 }
 \\ \nonumber
 &&\qquad\qquad\qquad\qquad\qquad\qquad\qquad\qquad\qquad
 \forall \tau\in  V^{-1,\alpha}(\partial\Omega)\,,
\end{eqnarray}
where $S_{\Omega,+}^t$ is the transpose map of $S_{\Omega,+}$.
\end{definition}
As shown in \cite[\S 13]{La24b},  $(V^{-1,\alpha}(\partial\Omega), \|\cdot\|_{  V^{-1,\alpha}(\partial\Omega)  })$ is a Banach space. By definition of the norm, $C^{0,\alpha}(\partial\Omega)$ is continuously embedded into $V^{-1,\alpha}(\partial\Omega)$. Moreover, we have the following statement  of \cite[Prop.~6.6]{La24c}  on the continuity of the normal derivative on $C^{0,\alpha}(\overline{\Omega})_\Delta$. 
\begin{proposition}\label{prop:ricodnu}
 Let   $\alpha\in ]0,1[$. Let $\Omega$ be a bounded open  subset of 
 ${\mathbb{R}}^{n}$ of class $C^{1,\alpha}$. Then the distributional normal derivative   $\partial_\nu$ is a continuous surjection of $C^{0,\alpha}(\overline{\Omega})_\Delta$ onto $V^{-1,\alpha}(\partial\Omega)$ and there exists $Z\in {\mathcal{L}}\left(V^{-1,\alpha}(\partial\Omega),C^{0,\alpha}(\overline{\Omega})_\Delta\right)$ such that
 \begin{equation}\label{prop:ricodnu1}
\partial_\nu Z[g]=g\qquad\forall g\in V^{-1,\alpha}(\partial\Omega)\,,
\end{equation}
\textit{i.e.}, $Z$ is a right inverse of  $\partial_\nu$. (See Lemma \ref{lem:c1alcof} for the topology of $C^{0,\alpha}_{
{\mathrm{loc}}	}(\overline{\Omega^-})_\Delta $).
\end{proposition}
We also mention that the condition of Definition \ref{defn:conoderdedu} admits a different formulation at least in case $u\in C^{0,\alpha}(\overline{\Omega})_\Delta$ (see \cite[Prop.~5.15]{La24c}).
\begin{proposition}\label{prop:node1adeq}
 Let   $\alpha\in ]0,1[$. Let $\Omega$ be a bounded open  subset of 
 ${\mathbb{R}}^{n}$ of class $C^{1,\alpha}$. Let $E_\Omega$ be a linear map from $C^{1,\alpha}(\partial\Omega)$ to $ C^{1,\alpha}(\overline{\Omega})_\Delta$ such that
 \begin{equation}\label{prop:node1adeq0}
 E_\Omega[f]_{|\partial\Omega}=f\qquad\forall f\in C^{1,\alpha}(\partial\Omega)\,.
 \end{equation}
 If $u\in C^{0,\alpha}(\overline{\Omega})_\Delta$, then the distributional  normal derivative $\partial_\nu u$
 of $u$  is characterized by the validity of the following equality
 \begin{eqnarray}\label{prop:node1adeq1}
 \lefteqn{
\langle \partial_\nu u ,v\rangle =\int_{\partial\Omega}u
\frac{\partial}{\partial\nu}E_\Omega[v]
\,d\sigma
}
\\ \nonumber
&&\qquad\qquad
+\langle E^\sharp_\Omega[\Delta u],E_\Omega[v]\rangle -\int_\Omega\Delta (E_\Omega[v]) u\,dx
\qquad\forall v\in C^{1,\alpha}(\partial\Omega)\,.
\end{eqnarray}
\end{proposition}
   We note that in equality (\ref{prop:node1adeq1}) we can use the `extension' operator $E_\Omega$ that we prefer and that accordingly equality (\ref{prop:node1adeq1}) is independent of the specific choice of $E_\Omega$. When we deal with problems for the Laplace operator, a good choice is
 $E_\Omega={\mathcal{G}}_{\Omega,d,+}$, so that the last term in the right hand side of (\ref{prop:node1adeq1}) disappears.  
 In particular, under the assumptions of Proposition \ref{prop:node1adeq}, an extension operator as $E_\Omega$ always exists.
 
 \section{A  distributional outward unit normal derivative for the functions of $C^{0,\alpha}_{
{\mathrm{loc}}	}(\overline{\Omega^-})_\Delta$}
\label{sec:doun}
 
 We now plan to define the normal derivative on the boundary the functions of $C^{0,\alpha}_{
{\mathrm{loc}}	}(\overline{\Omega^-})_\Delta$ with $\alpha\in]0,1[$ in case $\Omega$ is a bounded open subset of ${\mathbb{R}}^n$ of class $C^{1,\alpha}$ as in \cite[\S 3]{La25}.  To do so, we  choose $r\in]0,+\infty[$ such that $\overline{\Omega}\subseteq {\mathbb{B}}_n(0,r)$ and we observe that the linear operator from $C^{1,\alpha}(\partial\Omega)$ to $C^{1,\alpha}((\partial\Omega)\cup(\partial{\mathbb{B}}_n(0,r)))$ that is defined by the equality
\begin{equation}\label{eq:zerom}
\stackrel{o}{E}_{\partial\Omega,r}[v]\equiv\left\{
\begin{array}{ll}
 v(x) & \text{if}\ x\in \partial\Omega\,,
 \\
 0 & \text{if}\ x\in \partial{\mathbb{B}}_n(0,r)\,,
\end{array}
\right.\qquad\forall v\in C^{1,\alpha}(\partial\Omega) 
\end{equation}
is linear and continuous and that accordingly the transpose map $\stackrel{o}{E}_{\partial\Omega,r}^t$ is linear and continuous
\[
\text{from}\ 
\left(C^{1,\alpha}((\partial\Omega)\cup(\partial{\mathbb{B}}_n(0,r)))\right)'\ 
\text{to}\ \left(C^{1,\alpha}(\partial\Omega)\right)'\,.
\]
Then we introduce the following definition as in \cite[\S 3]{La25}.
\begin{definition}\label{defn:endedr}
 Let   $\alpha\in ]0,1[$. Let $\Omega$ be a bounded open  subset of 
 ${\mathbb{R}}^{n}$ of class $C^{1,\alpha}$. Let $r\in]0,+\infty[$ such that $\overline{\Omega}\subseteq {\mathbb{B}}_n(0,r)$. If $u\in C^{0,\alpha}_{
{\mathrm{loc}}	}(\overline{\Omega^-})_\Delta$, then we set
\begin{equation}\label{defn:endedr1}
 \partial_{\nu_{\Omega^-}}u
 \equiv 
 \stackrel{o}{E}_{\partial\Omega,r}^t\left[
\partial_{\nu_{{\mathbb{B}}_n(0,r)\setminus\overline{\Omega}}}u
\right]\,.
\end{equation}
 \end{definition}
 Under the assumptions of Definition \ref{defn:endedr},   $\partial_{\nu_{\Omega^-}}u$ is an element of $\left(C^{1,\alpha}(\partial\Omega)\right)'$. 
Under the assumptions of Definition \ref{defn:endedr}, there exists  a linear (extension) map $E_{{\mathbb{B}}_n(0,r)\setminus\overline{\Omega}}$ from $C^{1,\alpha}((\partial\Omega)\cup(\partial
{\mathbb{B}}_n(0,r)
))$ to $ C^{1,\alpha}(\overline{{\mathbb{B}}_n(0,r)}\setminus\Omega)_\Delta$ such that
 \begin{equation}\label{prop:node1adeq0r}
 E_{{\mathbb{B}}_n(0,r)\setminus\overline{\Omega}}[f]_{|(\partial\Omega)\cup(\partial
{\mathbb{B}}_n(0,r)
)}=f\qquad\forall f\in C^{1,\alpha}((\partial\Omega)\cup(\partial
{\mathbb{B}}_n(0,r)
))\,.
 \end{equation}
Then the definition of normal derivative on the boundary of 
${\mathbb{B}}_n(0,r)\setminus\overline{\Omega}$ implies that
\begin{eqnarray}\label{defn:endedr2}
\lefteqn{
\langle \partial_{\nu_{\Omega^-}}u,v\rangle 
=-\int_{\partial\Omega}u\frac{\partial}{\partial\nu_{\Omega}}E_{{\mathbb{B}}_n(0,r)\setminus\overline{\Omega}}[\stackrel{o}{E}_{\partial\Omega,r}[v]]\,d\sigma
}
\\ \nonumber
&&\qquad\qquad 
+\int_{\partial{\mathbb{B}}_n(0,r)}u\frac{\partial}{\partial\nu_{{\mathbb{B}}_n(0,r)}}E_{{\mathbb{B}}_n(0,r)\setminus\overline{\Omega}}[\stackrel{o}{E}_{\partial\Omega,r}[v]]\,d\sigma
\\ \nonumber
&&\qquad\qquad 
+\langle E^\sharp_{{\mathbb{B}}_n(0,r)\setminus\overline{\Omega}}[\Delta u],E_{{\mathbb{B}}_n(0,r)\setminus\overline{\Omega}}[\stackrel{o}{E}_{\partial\Omega,r}[v]]\rangle 
\\ \nonumber
&&\qquad\qquad 
-\int_{{\mathbb{B}}_n(0,r)\setminus\overline{\Omega}}\Delta (E_{{\mathbb{B}}_n(0,r)\setminus\overline{\Omega}}[\stackrel{o}{E}_{\partial\Omega,r}[v]]) u\,dx
\qquad\forall v\in C^{1,\alpha}(\partial\Omega)\,.
\end{eqnarray}
 (cf.~(\ref{prop:node1adeq1})).   
As shown in \cite[\S 3]{La25},  Definition  \ref{defn:endedr} is independent of the choice of $r\in]0,+\infty[$ such that $\overline{\Omega}\subseteq {\mathbb{B}}_n(0,r)$.  

\begin{remark}\label{rem:conoderdeducl-}
 Let   $\alpha\in ]0,1[$. Let $\Omega$ be a bounded open  subset of 
 ${\mathbb{R}}^{n}$ of class $C^{1,\alpha}$, $u\in C^{1,\alpha}_{
{\mathrm{loc}}	}(\overline{\Omega^-})$. Then equality (\ref{defn:endedr1}) of  the  
 of Definition \ref{defn:endedr} of normal derivative and equality (\ref{lem:conoderdeducl1}) imply that
 \begin{equation}\label{rem:conoderdeducl-1}
 \langle \partial_{\nu_{\Omega^-}} u ,v\rangle \equiv\int_{\partial\Omega}\frac{\partial u}{\partial\nu_{\Omega^-}}v\,d\sigma
  \quad\forall v\in C^{1,\alpha}(\partial\Omega)\,,
\end{equation}
where $\frac{\partial u}{\partial\nu_{\Omega^-}}$ in the right hand side denotes the classical $\nu_{\Omega^-}$-normal derivative of $u$ on $\partial\Omega$.  In the sequel, we use the classical symbol $\frac{\partial u}{\partial\nu_{\Omega^-}}$ also for  $\partial_{\nu_{\Omega^-}} u$ when no ambiguity can arise.\par 
\end{remark}

Next we note that if  $u\in C^{0,\alpha}_{{\mathrm{loc}}}(\overline{\Omega^-})$ and $u$ is both harmonic in $\Omega^-$ and harmonic at infinity, then $u\in C^{0,\alpha}_{
{\mathrm{loc}}	}(\overline{\Omega^-})_\Delta$. As shown in \cite[\S 3]{La25}, 
  the distribution $\partial_\nu u $ of Definition \ref{defn:endedr} coincides with the normal derivative that has been introduced in \cite[Defn.~6.4]{La24b} for harmonic functions.   By \cite[\S 3]{La25}, we have the following continuity statement for the distributional normal derivative $\partial_{\nu_{\Omega^-}}$ of
Definition \ref{defn:endedr}.

\begin{proposition}\label{prop:recodnu}
 Let   $\alpha\in ]0,1[$. Let $\Omega$ be a bounded open  subset of 
 ${\mathbb{R}}^{n}$ of class $C^{1,\alpha}$.  Then the distributional normal derivative   $\partial_{\nu_{\Omega^-}}$ is a continuous surjection of $C^{0,\alpha}_{
{\mathrm{loc}}	}(\overline{\Omega^-})_\Delta$ onto $V^{-1,\alpha}(\partial\Omega)$ and there exists a linear and continuous  operator   $Z_-$ from  $ V^{-1,\alpha}(\partial\Omega)$ to $C^{0,\alpha}_{
{\mathrm{loc}}	}(\overline{\Omega^-})_\Delta $ such that
 \begin{equation}\label{prop:recodnu1}
\partial_{\nu_{\Omega^-}} Z_-[g]=g\qquad\forall g\in V^{-1,\alpha}(\partial\Omega)\,,
\end{equation}
\textit{i.e.}, $Z_-$ is a right inverse of  $\partial_{\nu_{\Omega^-}}$. (See Lemma \ref{lem:c1alcof} for the topology of $C^{0,\alpha}_{
{\mathrm{loc}}	}(\overline{\Omega^-})_\Delta $).	 
\end{proposition}

\section{Preliminaries on the fundamental solution of the Helmholtz equation}\label{sec:prefusoh}
Next we turn to introduce the fundamental solution of $\Delta +k^2$ when  $k\in{\mathbb{C}}\setminus]-\infty,0]$. In the sequel, $\arg$ and $\log$ denote the principal branch of the argument and of the logarithm in $ {\mathbb{C}}\setminus]-\infty,0]$, respectively. Then we have
\[
\arg (z)={\mathrm{Im}}\,  \log (z)\in]-\pi,\pi[\qquad\forall 
z\in {\mathbb{C}}\setminus]-\infty,0]\,.
\]
Then we set 
\begin{equation}
\label{sh1}
J^{\sharp}_{\nu}(z)\equiv
\sum_{j=0}^{\infty}\frac{
(-1)^{j}
z^{j}(1/2)^{2j}(1/2)^{\nu}
}{\Gamma(j+1)\Gamma(j+\nu+1)}\qquad\forall z\in {\mathbb{C}}\,,
\end{equation}
 for all $\nu\in{\mathbb{C}}\setminus\{-j:\,j\in {\mathbb{N}}\setminus\{0 \}
\}$. 		  	  Here $(1/2)^{\nu}=e^{\nu \log (1/2)}$.  	 
As is well known, if $\nu\in{\mathbb{C}}\setminus\{-j:\,j\in {\mathbb{N}}\setminus\{0 \}
\}$ then the function $J^{\sharp}_{\nu}(\cdot)$ is entire and 
\begin{equation}
\label{nens0}
J^{\sharp}_{\nu}(z^{2})\equiv e^{-\nu\log z}J_{\nu}(z)
\qquad\forall  z\in {\mathbb{C}}\setminus]-\infty,0] \,,		 
\end{equation}
where $J_{\nu}(\cdot)$ is the Bessel function of the first kind of index $\nu$ (cf.~\textit{e.g.}, Lebedev~\cite[Ch.~1, \S 5.3]{Le72}.)  One could also consider case $\nu\in -{\mathbb{N}}$, but we do not need such a case in this paper.		 		If $\nu\in {\mathbb{N}}$, we set
\begin{eqnarray}
\label{nens1}
\lefteqn{
N^{\sharp}_{\nu}(z)\equiv 
-\frac{2^{\nu}}{\pi}\sum_{0\leq j\leq \nu-1}
\frac{(\nu-j-1)!}{j!}z^{j}(1/2)^{2j }
}
\\
\nonumber
&&
\qquad
-\frac{z^{\nu}}{\pi}
\sum_{j=0}^{\infty}\frac{(-1)^{j}z^{j}(1/2)^{2j}(1/2)^{\nu}}{j!(\nu+j)!}
\left(
2\sum_{0<l\leq j} \frac{1}{l}+\sum_{j<l\leq j+\nu}\frac{1}{l}
\right)\quad\forall z\in {\mathbb{C}}\,.
\end{eqnarray}
As one can easily see, the    $N^{\sharp}_{\nu}(\cdot)$ is  an entire holomorphic function of the variable $z\in {\mathbb{C}}$ and 
\begin{equation}
\label{nenent1}
N_{\nu}(z)=\frac{2}{\pi}
(\log (z)-\log 2+\gamma)J_{\nu}(z)+z^{-\nu}N^{\sharp}_{\nu}(z^{2})
\qquad\forall z \in {\mathbb{C}}\setminus]-\infty,0]\,,
\end{equation}
where $\gamma$ is the Euler-Mascheroni constant, and where
$N_{\nu}(\cdot)$ is the Neumann function of index $\nu$, also known as Bessel function of second kind and index $\nu$ (cf.~\textit{e.g.}, ~Lebedev~\cite[Ch.~1, \S 5.5]{Le72}.) Let   $k\in{\mathbb{C}}\setminus]-\infty,0]$, $n\in {\mathbb{N}}\setminus\{0,1\}$,  $a_{n}\in {\mathbb{C}}$.  Then we set
\begin{equation}
\label{bn}
b_{n}\equiv
\left\{
\begin{array}{ll}
\pi^{1-(n/2)}2^{-1-(n/2)} & {\text{if}}\ n\  {\text{is\ even}}\,,
\\
(-1)^{\frac{n-1}{2}}\pi^{1-(n/2)}2^{-1-(n/2)} & {\text{if}}\ n\ {\text{is\ odd}}\,,
\end{array}
\right.
\end{equation}
and
\begin{equation}
\label{urassl2}
\tilde{S}_{ n,  k,a_{n} } (x)=\left\{
\begin{array}{ll}
k^{ n-2}
\biggl\{
a_{n}+\frac{2b_{n}}{\pi} (
 \log \,k
-\log 2+\gamma)+
\frac{2b_{n}}{\pi} \log|x|
\biggr\}      &
\\
\quad  \times J_{\frac{n-2}{2}}^{\sharp}(k^2 |x|^{2})
+b_{n}|x|^{2-n}N_{\frac{n-2}{2}}^{\sharp}(k^2|x|^{2})  &
\\
\quad\qquad\qquad\qquad\quad\ 
 \qquad\qquad\qquad \qquad \qquad   {\text{if}}\   n \ {\text{is\ even}}, &
\\
a_{n}k^{ n-2 }J_{\frac{n-2}{2}}^{\sharp}(k^2|x|^{2})+b_{n}|x|^{2-n}J_{-\frac{n-2}{2}}^{\sharp}(k^2|x|^{2})&
\\
\quad\qquad\qquad\qquad\quad\ \
 \qquad\qquad\qquad \qquad \qquad  {\text{if}}\  n\ {\text{is\ odd}}, &
\end{array}
\right.
\end{equation}
for all $x\in{\mathbb{R}}^{n}\setminus\{0\}$. As it is known and can be easily verified, the family $\{\tilde{S}_{  n, k,a_{n} } \}_{a_{n}\in {\mathbb{C}} }$ coincides with the family of all radial fundamental solutions of $\Delta+k^2$.

Now let $k\in {\mathbb{C}}\setminus]-\infty,0]$, ${\mathrm{Im}}\, k\geq 0$. As well known in scattering theory, a function $u$ of class $C^1$ in the complement of a compact subset of ${\mathbb{R}}^{n}$ satisfies the outgoing $k$-radiation condition  provided that
\begin{equation}
\label{rad1}
\lim_{x\to\infty}|x|^{\frac{n-1}{2}}(Du(x)\frac{x}{|x|}-iku(x))=0\,.
\end{equation}
Now by writing the fundamental solution of (\ref{urassl2}) 
in terms of the Hankel functions, and by exploiting the asymptotic behavior at infinity of the Hankel functions, one finds classically that the fundamental solution of (\ref{urassl2}) satisfies the outgoing $k$-radiation condition if and only if
\begin{equation}
\label{raso1aaa}
a_{n}\equiv\left\{
\begin{array}{ll}
-ib_{n}
& {\mathrm{if}}\ n\ {\text{is\ even}} 
\\
-e^{-i\frac{n-2}{2}\pi}b_{n}
& {\mathrm{if}}\ n\ {\text{is\ odd}}
\end{array}
\right\}=-i\pi^{1-(n/2)}2^{-1-(n/2)}\,.
\end{equation}
Then we introduce the following definition.
\begin{definition}
\label{defn:raso}
Let $k\in {\mathbb{C}}\setminus]-\infty,0]$, ${\mathrm{Im}}\, k\geq 0$. We denote by $\tilde{S}_{n,k;r}$ the function from  ${\mathbb{R}}^{n}\setminus\{0\}$ to ${\mathbb{C}}$ defined by 
\[
\tilde{S}_{n,k;r}(x)\equiv
\tilde{S}_{n,  k,a_{n} } (x)
\qquad\forall x\in {\mathbb{R}}^{n}\setminus\{0\}\,,
\]
with   $a_{n}$   as in (\ref{raso1aaa}) (cf. (\ref{urassl2}).)
\end{definition}
As we have said above, if $k\in {\mathbb{C}}\setminus]-\infty,0]$ and ${\mathrm{Im}}\, k\geq 0$, then $\tilde{S}_{n,k;r}$
satisfies the outgoing $k$-radiation condition   and $\tilde{S}_{n,k;r}$ is the fundamental solution that is classically used in scattering theory. The subscript $r$ stands for `radiation'.

In particular, for $n=3$  we have
$\tilde{S}_{n,k;r}(x)=-\frac{1}{4\pi |x|}e^{ik|x|}$ for all $x\in {\mathbb{R}}^{n}\setminus\{0\}$. Then by exploiting the definition of the first Hankel function
\begin{equation}\label{eq:hank1}
H^{(1)}_{ \frac{n-2}{2}}(z)\equiv J_{ \frac{n-2}{2}}(z)+i N_{ \frac{n-2}{2}}(z)\qquad\forall z\in {\mathbb{C}}\setminus]-\infty,0]\,,
\end{equation}
we have
\begin{equation}\label{eq:sofura}
\tilde{S}_{n,k;r}(x)=C_n|x|^{-\frac{n-2}{2}}H^{(1)}_{ \frac{n-2}{2}}(k|x|)
\qquad\forall x\in {\mathbb{R}}^n\setminus\{0\}\,,
\end{equation}
where 
\begin{equation}\label{eq:sofura1}
C_{n}\equiv\left\{
\begin{array}{ll}
\frac{	k^{ \frac{n-2}{2} } }{i4(2\pi)^{\frac{n-2}{2}}}
& {\text{if}}\ n\ {\text{is\ even}}\,,
\\
\frac{e^{  \frac{n-2}{2}  \log (k)}}{i4(2\pi)^{\frac{n-2}{2}}}
& {\text{if}}\ n\ {\text{is\ odd}}\,,
\end{array}
\right.
\end{equation}
 (see also Mitrea \cite[(7.6.11), p.~266]{Mit18}. Mitrea, Mitrea and Mitrea \cite[(6.1.7), p.~882]{MitMitMit23}).  Then known properties of the first Hankel function imply the validity of the following statement. 
\begin{proposition}\label{prop:furasoin}
 Let $k\in {\mathbb{C}}\setminus]-\infty,0]$, ${\mathrm{Im}}\, k\geq 0$. Let
 \begin{equation}\label{prop:furasoina}
\eta_n(\rho)\equiv C_n\rho^{-\frac{n-2}{2}}H^{(1)}_{ \frac{n-2}{2}}(k\rho)
\qquad\forall \rho\in]0,+\infty[\,.
\end{equation}
Then the following statements hold true.
\begin{enumerate}
\item[(i)] 
\begin{equation}\label{prop:furasoin0}
 a_0\equiv\sup_{\rho\in[1,+\infty[}\rho^{\frac{n-1}{2}}e^{{\mathrm{Im}}\,(k)\rho}\eta_n(\rho) <+\infty\,.
 \end{equation}

\item[(ii)] $\eta_{n}'(\rho)=(-k)\rho^{-\frac{n-2}{2}}C_{n}H^{(1)}_{\frac{n-2}{2}+1}(k\rho)
$ for all $\rho\in]0,+\infty[$ and
 \begin{equation}\label{prop:furasoin1}
a_1\equiv\sup_{\rho\in[1,+\infty[}\rho^{\frac{n-1}{2}}e^{{\mathrm{Im}}\,(k)\rho}|\eta_n'(\rho) |<+\infty\,.
\end{equation}

\item[(iii)] 
 \[
\rho^{\frac{n-1}{2}}(\eta'_{n}(\rho)- ik\eta_{n}(\rho))
=-k\rho^{1/2} 
C_{n}\biggl(
H^{(1)}_{\frac{n-2}{2}+1}(k\rho)+ iH^{(1)}_{\frac{n-2}{2}}(k\rho)
\biggr)
\]
for all $\rho\in]0,+\infty[$ and
\begin{equation}
 \label{prop:furasoin2}
 a_2\equiv\sup_{\rho\in[1,+\infty[}\rho^{\frac{n+1}{2}}e^{{\mathrm{Im}}\,(k)\rho}
 |\eta_n'(\rho)-ik\eta_n(\rho)|<+\infty\,.
\end{equation}
\item[(iv)] \[
\eta_n''(\rho)=
-k\rho^{-\frac{n-2}{2}-1} 
 C_{n}
H^{(1)}_{\frac{n-2}{2}+1}(k\rho)
+k^{2}\rho^{-\frac{n-2}{2}}  
 C_{n}
H^{(1)}_{\frac{n-2}{2}+2}(k\rho)
\ \ \forall \rho\in]0,+\infty[
\]
and
\begin{equation}
\label{prop:furasoin3}
 a_3\equiv\sup_{\rho\in[1,+\infty[}\rho^{\frac{n-1}{2}}e^{{\mathrm{Im}}\,(k)\rho}|\eta_n''(\rho) |<+\infty\,.
\end{equation}

\item[(v)] \begin{eqnarray*}
\lefteqn{
\rho^{\frac{n-1}{2}}(\eta''_{n}(\rho)- ik\eta'_{n}(\rho))
}
\\  \nonumber
&&
\ 
=\rho^{1/2}(-k)^{2}
C_{n}\biggl(-(k\rho)^{-1}H^{(1)}_{\frac{n-2}{2}+1}(k\rho)
+
H^{(1)}_{\frac{n-2}{2}+2}(k\rho)
+ i
H^{(1)}_{\frac{n-2}{2}+1}(k\rho)
\biggr)
\end{eqnarray*}
for all $\rho\in]0,+\infty[$ and
\begin{equation}
\label{prop:furasoin4}
 a_4\equiv\sup_{\rho\in[1,+\infty[}\rho^{\frac{n+1}{2}}e^{{\mathrm{Im}}\,(k)\rho}|\eta_n''(\rho)-ik\eta_n'(\rho)|<+\infty\,.
\end{equation}
\item[(vi)]
\[
\eta_{n}'''(\rho)=3k^2\rho^{-\frac{n-2}{2}-1} 
C_{n}H^{(1)}_{\frac{n-2}{2}+2}(k\rho)
-k^3\rho^{-\frac{n-2}{2}} 
C_{n}H^{(1)}_{\frac{n-2}{2}+3}(k\rho)
\]
for all $\rho\in]0,+\infty[$ and
\begin{equation}
\label{prop:furasoin5}
 a_5\equiv\sup_{\rho\in[1,+\infty[}\rho^{\frac{n-1}{2}}e^{{\mathrm{Im}}\,(k)\rho}|\eta_n'''(\rho) |<+\infty\,.
\end{equation}
\item[(vii)]
\begin{eqnarray*}
\lefteqn{
\rho^{\frac{n-1}{2}}(\eta'''_{n}(\rho)- ik\eta''_{n}(\rho))
}
\\ \nonumber
&&\qquad 
=\rho^{-1/2} 
\biggl[3k^2 C_{n}
H^{(1)}_{\frac{n-2}{2}+2}(k\rho)
 + ik^2
 C_{n}
H^{(1)}_{\frac{n-2}{2}+1}(k\rho)
 \biggr]
\\ \nonumber
&&\qquad\quad
-\rho^{1/2} k^3\biggl[
C_n\biggl(H^{(1)}_{\frac{n-2}{2}+3}(k\rho)+ iH^{(1)}_{\frac{n-2}{2}+2}(k\rho)\biggr)
\biggr]
\end{eqnarray*}
for all $\rho\in]0,+\infty[$ and
\begin{equation}
\label{prop:furasoin6}
 a_6\equiv\sup_{\rho\in[1,+\infty[}\rho^{\frac{n+1}{2}}e^{{\mathrm{Im}}\,(k)\rho}|\eta_n'''(\rho)-ik\eta_n''(\rho)|<+\infty\,.
\end{equation}
\end{enumerate}
 \end{proposition} 
 {\bf Proof.} All formulas for the derivatives of $\eta_n$ follow by the classical differentiation rules and by exploiting the following formula for the derivative  of the first Hankel function 
 \[
 \frac{d}{d\rho}\left(  \rho^{-\nu}H^{(1)}_{\nu}(k\rho)\right)
 =-k\rho^{-\nu}H^{(1)}_{\nu+1}(k\rho)\qquad\forall \rho\in]0,+\infty[\,,
 \] 
  for all $\nu\in\left\{\frac{n-2}{2},\frac{n-2}{2}+1, \frac{n-2}{2}+2\right\}$
 (cf.~\textit{e.g.}, 
 Lebedev~\cite[(5.6.3), p.~108]{Le72}). 
 
 The inequalities for the derivatives of $\eta_n$ of statements 
 (i), (ii), (iv), (vi) follow  by exploiting equality (\ref{prop:furasoina}),  
 the equalities for the derivatives of  statements 
  (ii), (iv), (vi) and the following inequality
  \[
  \sup_{\rho\in[1,+\infty[}\rho^{1/2}e^{{\mathrm{Im}}\,(k)\rho}\left|H^{(1)}_{\nu}(k\rho)\right|<+\infty\,,
  \]
  for all $\nu\in\left\{\frac{n-2}{2},\frac{n-2}{2}+1, \frac{n-2}{2}+2\right\}$, 
 which follows by the expansion in Lebedev~\cite[(5.11.4), p.~122]{Le72}. 
 
 In order to prove the inequalities of statements
 (iii), (v), (vii), we exploit equality (\ref{prop:furasoina}), the equalities for the derivatives of $\eta_n$ of statements (ii), (iv), (vi) and the expansion in Lebedev~\cite[(5.11.4), p.~122]{Le72}. Here, we note that  the first order term of $H^{(1)}_{\nu+1}(k\rho)$ cancels out with the first order term  of $iH^{(1)}_{\nu}(k\rho)$, so that
 \[
 \sup_{\rho\in[1,+\infty[}\rho e^{{\mathrm{Im}}\,(k)\rho}\rho^{1/2}\left|H^{(1)}_{\nu+1}(k\rho) + iH^{(1)}_{\nu}(k\rho))\right|<+\infty
 \]
     for all $\nu\in\left\{\frac{n-2}{2},\frac{n-2}{2}+1, \frac{n-2}{2}+2\right\}$. \hfill  $\Box$ 

\vspace{\baselineskip}

 \section{Acoustic layer potentials with densities in the space $V^{-1,\alpha}(\partial\Omega)$}
 \label{sec:acsilah-1a}
 
 Let  $\alpha\in]0,1]$. Let  $\Omega$ be a bounded open subset of ${\mathbb{R}}^{n}$ of class $C^{1,\alpha}$. Let  $r_{|\partial\Omega}$  be the restriction map  from ${\mathcal{D}}({\mathbb{R}}^n)$ to $C^{1,\alpha}(\partial\Omega)$.  Let $\lambda\in {\mathbb{C}}$. If $S_{n,\lambda} $ is a fundamental solution 
of the operator $\Delta+\lambda$ and $\mu\in (C^{1,\alpha}(\partial\Omega))'$, then the  (distributional) single layer  potential relative to $S_{n,\lambda} $ and $\mu$ is the distribution
\[
v_\Omega[S_{n,\lambda} ,\mu]=(r_{|\partial\Omega}^t\mu)\ast S_{n,\lambda}  \in {\mathcal{D}}'({\mathbb{R}}^n)\,.
\]
 and we  set
\begin{eqnarray}\label{eq:dsila}
v_\Omega^+[S_{n,\lambda} ,\mu]&\equiv&\left((r_{|\partial\Omega}^t\mu)\ast S_{n,\lambda}  \right)_{|\Omega}
\qquad\text{in}\ \Omega\,,
\\ \nonumber
v_\Omega^-[S_{n,\lambda} ,\mu] &\equiv&
\left((r_{|\partial\Omega}^t\mu)\ast S_{n,\lambda}  \right)_{|\Omega^-}
\qquad\text{in}\ \Omega^-\,.
\end{eqnarray}
It is also known that the restriction of $v_\Omega[S_{n,\lambda} ,\mu]$ to ${\mathbb{R}}^n\setminus\partial\Omega$ equals the (distribution that is associated to) the function
\[
 \langle (r_{|\partial\Omega}^t\mu)(y),S_{n,\lambda} (\cdot-y)\rangle \,.
\]
In the cases in which both $v_\Omega^+[S_{n,\lambda} ,\mu]$ and $v_\Omega^-[S_{n,\lambda} ,\mu]$ admit a continuous extension to $\overline{\Omega}$ and to $\overline{\Omega^-}$, respectively, we still use the symbols $v_\Omega^+[S_{n,\lambda} ,\mu]$ and $v_\Omega^-[S_{n,\lambda} ,\mu]$ for the continuous extensions and if the values of $v_\Omega^\pm[S_{n,\lambda} ,\mu](x)$ coincide for each $x\in\partial\Omega$, then we set
\[
V_\Omega^+[S_{n,\lambda} ,\mu](x)\equiv v_\Omega^+[S_{n,\lambda} ,\mu](x)=v_\Omega^-[S_{n,\lambda} ,\mu]^-(x)\qquad\forall x\in\partial\Omega\,.
\]
If $\mu$ is continuous, then it is known that $v_\Omega[S_{n,\lambda} ,\mu ]$ is continuous in ${\mathbb{R}}^n$. Indeed,  $\partial\Omega$ is upper $(n-1)$-Ahlfors regular with respect to ${\mathbb{R}}^n$ and $S_{n,\lambda}$ has a weak singularity (cf.~\cite[Prop.~6.5]{La24}, \cite[Prop.~4.3]{La22b}, \cite[Lem.~4.2 (i)]{DoLa17}). If $\lambda=0$, \textit{i.e.} if $\Delta+\lambda $ is the Laplace operator, and if we take the fundamental solution $S_n$ of $\Delta$ that is delivered by the formula
 \[
S_{n}(\xi)\equiv
\left\{
\begin{array}{lll}
\frac{1}{s_{n}}\ln  |\xi| \qquad &   \forall \xi\in 
{\mathbb{R}}^{n}\setminus\{0\},\quad & {\mathrm{if}}\ n=2\,,
\\
\frac{1}{(2-n)s_{n}}|\xi|^{2-n}\qquad &   \forall \xi\in 
{\mathbb{R}}^{n}\setminus\{0\},\quad & {\mathrm{if}}\ n>2\,,
\end{array}
\right.
\]
where $s_{n}$ denotes the $(n-1)$ dimensional measure of 
$\partial{\mathbb{B}}_{n}(0,1)$, then  we set
\begin{equation}\label{eq:vosn}
v_\Omega[\mu]=v_\Omega[S_n,\mu]\,,\qquad
V_\Omega[\mu]=V_\Omega[S_n,\mu]
\,.
\end{equation}
  Next we introduce the following (classical) technical statement on the  double layer potential (for a proof we refer, \textit{e.g.}, to \cite[\S 5]{La25a}).
\begin{theorem}\label{thm:dlay}
 Let  $\alpha\in]0,1[$. Let $\Omega$ be a bounded open subset of  ${\mathbb{R}}^n$ of class $C^{1,\alpha}$. Let $\lambda\in {\mathbb{C}}$.  Let $S_{n,\lambda} $ be a fundamental solution 
of the operator $\Delta+\lambda$. If $\mu\in C^{0,\alpha}(\partial\Omega)$, and
\begin{equation}\label{thm:dlay1}
w_\Omega[S_{n,\lambda} ,\mu](x)\equiv\int_{\partial\Omega}\frac{\partial}{\partial\nu_{\Omega,y}}\left(S_{n,\lambda} (x-y)\right)\mu(y)\,d\sigma_y\qquad\forall x\in {\mathbb{R}}^n\,,
\end{equation}
where
\[
\frac{\partial}{\partial \nu_{\Omega_y} }
\left(S_{n,\lambda}(x-y)\right)\equiv
  -DS_{n,\lambda}(x-y) \nu_{\Omega}(y) \qquad\forall (x,y)\in\mathbb{R}^n\times\partial\Omega\,, x\neq y\,,
 \]
then the restriction 
$w_\Omega[S_{n,\lambda} ,\mu]_{|\Omega}$ can be extended uniquely to a function  
$w^{+}_\Omega [S_{n,\lambda} ,\mu]$ of class $ C^{0,\alpha}(\overline{\Omega})$ and $w_\Omega[S_{n,\lambda} ,\mu]_{|{\mathbb{R}}^n\setminus\overline{\Omega}}$ can be extended uniquely to a function  
$w^{-}_\Omega [S_{n,\lambda} ,\mu]$ of class  $
C^{0,\alpha}_{ {\mathrm{loc}} }(\overline{\Omega^{-}})$.\par 

Moreover, the map 
from   $C^{0,\alpha}(\partial\Omega)$ to    $C^{0,\alpha}(\overline{\Omega})$
 which takes $\mu$ to $
w^{+}_\Omega [S_{n,\lambda} ,\mu]$ is  continuous and the map from     $C^{0,\alpha}(\partial\Omega)$ to   $C^{0,\alpha}(\overline{{\mathbb{B}}_{n}(0,r)}\setminus \Omega) $ which takes $\mu$ to 
$w^{-}_\Omega [S_{n,\lambda} ,\mu]_{|\overline{{\mathbb{B}}_{n}(0,r)}\setminus \Omega}$ is continuous for all $r\in]0,+\infty[$ such that $\overline{\Omega}\subseteq {\mathbb{B}}_{n}(0,r)$  and we have the following jump relation
\begin{equation}\label{thm:dlay1a}
w^{\pm}_\Omega [S_{n,\lambda},\mu](x)
=\pm\frac{1}{2}\mu(x)+w_\Omega[S_{n,\lambda},\mu](x)
\qquad\forall x\in\partial\Omega\,.
\end{equation}
\end{theorem}
 We also set
\[
W_\Omega[S_{n,\lambda} ,\mu](x)\equiv w_\Omega[S_{n,\lambda} ,\mu](x)\qquad\forall x\in\partial\Omega\,.
\]
If $\lambda=0$, \textit{i.e.} if $\Delta+\lambda $ is the Laplace operator, and if we take the fundamental solution $S_n$ of $\Delta$ , then  we set
\begin{equation}\label{eq:wosn}
w_\Omega[\mu]\equiv w_\Omega[S_n,\mu]\,,\qquad W_\Omega[\mu]\equiv W_\Omega[S_n,\mu]\,.
\end{equation}
 Next we state the following continuity statement for the acoustic single layer potential. For a proof, we refer to \cite[\S 5]{La25a}.
\begin{theorem}\label{thm:slhco-1a} 
Let $\alpha\in]0,1[$. Let $\Omega$ be a bounded open subset of ${\mathbb{R}}^n$ of class $C^{1,\alpha}$. Let $\lambda\in {\mathbb{C}}$.  Let $S_{n,\lambda} $ be a fundamental solution 
of the operator $\Delta+\lambda$. Let $r\in]0,+\infty[$ be such that $\overline{\Omega}\subseteq 
{\mathbb{B}}(0,r)$.

 If $\mu\in V^{-1,\alpha}(\partial\Omega)$, then $v_\Omega^+[S_{n,\lambda} ,\mu]\in C^{0,\alpha}(\overline{\Omega})_\Delta$ and $v_\Omega^-[S_{n,\lambda} ,\mu]\in C^{0,\alpha}_{{\mathrm{loc}}	}(\overline{\Omega^-})_\Delta$. Moreover, the map from $V^{-1,\alpha}(\partial\Omega)$ to $C^{0,\alpha}(\overline{\Omega})_\Delta$ that takes $\mu$ to
 $v_\Omega^+[S_{n,\lambda} ,\mu]$ is continuous 
 and
 the map from $V^{-1,\alpha}(\partial\Omega)$ to $C^{0,\alpha}(\overline{{\mathbb{B}}(0,r)}\setminus\Omega)_\Delta$ that takes $\mu$ to
 $v_\Omega^-[S_{n,\lambda} ,\mu]_{|\overline{{\mathbb{B}}(0,r)}\setminus\Omega}$ is continuous.
 \end{theorem}

\section{Acoustic layer potentials and the radiation condition}\label{sec:laporaco}
 
As is well-known,    the kernels of the simple and double layer potentials associated to the radial fundamental solution satisfy 
the outgoing $(k)$-radiation condition.  For the convenience of the reader, we include a statement in the form we need and a proof.
\begin{proposition}
\label{prop:dora}
Let  $k\in {\mathbb{C}}\setminus]-\infty,0]$, ${\mathrm{Im}}\,k\geq 0$. Let $K$ be a compact subset of ${\mathbb{R}}^{n}$. Let $v$ be a  function of $K$ to $\partial{\mathbb{B}}_{n}(0,1)$.  Then the following  inequalities hold
\begin{eqnarray}
\label{prop:dora0}
\lefteqn{
c_{0,K}\equiv \sup_{\zeta\in {\mathbb{N}}^n, |\zeta|\leq 2}\,
\sup_{|x|\geq 1+2\sup_{y\in K}|y|}\,
\sup_{y\in K}
\biggl\{\biggr.
e^{({\mathrm{Im}}\,k)|x|}|x||x|^{\frac{n-1}{2}}
}
\\ \nonumber
&&\qquad 
\times\biggl|
\frac{x}{|x|} D_{x} D^\zeta_{x}\tilde{S}_{n,k;r} (x-y)-(ik)D^\zeta_{x}\tilde{S}_{n,k;r}(x-y)
\biggr|\biggl.\biggr\}<+\infty\,,
\\
\label{prop:dora1}
\lefteqn{
c_{1,K}\equiv  \sup_{|x|\geq 1+2\sup_{y\in K}|y|}\,\sup_{y\in K}\biggl\{\biggr.e^{({\mathrm{Im}}\,k)|x|}|x||x|^{\frac{n-1}{2}}
}
\\ \nonumber
&&\qquad 
\times\biggl|
\frac{x}{|x|} D_{x}\biggl(
\frac{\partial}{\partial v(y)}
( \tilde{S}_{n,k;r} (x-y))
\biggr)
-(ik)
\frac{\partial}{\partial v(y)}
(\tilde{S}_{n,k;r}(x-y))
\biggr|\biggl.\biggr\}<+\infty\,.
\end{eqnarray}
Here
$\frac{\partial}{\partial v(y)}
(\tilde{S}_{n,k;r}(x-y))\equiv - v(y) \cdot D\tilde{S}_{n,k;r} (x-y)$.
\end{proposition}
{\bf Proof.} We first prove inequality (\ref{prop:dora0}). It clearly suffices to estimate
the argument of the supremum in (\ref{prop:dora0}) separately when
 $|\zeta|=0$, $|\zeta|=1$ and $|\zeta|=2$. In case $\zeta=0$,
we note that
\begin{eqnarray}
\label{prop:dora2}
\lefteqn{
\biggl|
\frac{x}{|x|} D_{x}\tilde{S}_{n,k;r}(x-y)-(ik)\tilde{S}_{n,k;r}(x-y)
\biggr|
}
\\ \nonumber
&& \leq
\biggl|\frac{x}{|x|}\cdot \frac{x-y}{|x-y|} \eta_n'(|x-y|)-ik\eta_n(|x-y|)
\biggr|
\\ \nonumber
&& \leq
\biggl|\frac{x}{|x|}\cdot \frac{x-y}{|x-y|} 
[\eta_n'(|x-y|)-ik\eta_n(|x-y|)]
\biggr|
\\ \nonumber
&& \quad
+
\biggl|ik\eta_n(|x-y|)
\biggl[
\frac{x}{|x|}\cdot \frac{x-y}{|x-y|} -1
\biggr]
\biggr|
\\ \nonumber
&& \leq
|\eta_n'(|x-y|)-ik\eta_n(|x-y|)|
+
|k||\eta_n(|x-y|)|\left|\frac{x}{|x|}-\frac{x-y}{|x-y|}\right|\,,
\end{eqnarray}
for all $y\in K$ and $x\in {\mathbb{R}}^{n}\setminus{\mathbb{B}}_{n}(0,1+\sup_{y\in K}|y|)$.
In order to estimate the right hand side of (\ref{prop:dora2}), we introduce the following two elementary inequalities.
\begin{eqnarray}
\label{prop:dora3}
\lefteqn{
\sup_{y\in K}\frac{|x|}{|x-y|}\leq
1+\sup_{y\in K}\frac{|y|}{|x-y|}
}
\\ \nonumber
&&\qquad\leq
1+\sup_{y\in K}\frac{|y|}{||x|-|y||}
\leq
1+\sup_{y\in K}|y|\qquad
\forall x\in {\mathbb{R}}^{n}\setminus{\mathbb{B}}_{n}(0,1+\sup_{y\in K}|y|)
\,,
\end{eqnarray}
and
\begin{eqnarray}
\label{prop:dora4}
\lefteqn{
\left|\frac{x}{|x|}-\frac{x-y}{|x-y|}\right|
=
\frac{|\,(x-y)|x|-x|x-y|\,|}{|x||x-y|}
}
\\ \nonumber
&&\qquad\leq
 \frac{|\,x|x|-x|x-y|\,|}{|x||x-y|}
 +
  \frac{|-y |x|\,|}{|x||x-y|}
\\ \nonumber
&& \qquad\leq
\frac{|\,|x|-|x-y|\,|}{|x-y|}
+
\frac{|y|}{|x-y|}  
\\ \nonumber
&&\qquad \leq
2\frac{|y|}{|x-y|}\leq
2\frac{|y|}{|x|-\sup_{y\in K}|y|}=2\frac{|y|}{2^{-1}|x|+2^{-1}|x|-\sup_{y\in K}|y|}
\\ \nonumber
&&\qquad\leq
4\frac{\sup_{y\in K}|y|}{|x|}
\qquad
\forall x\in {\mathbb{R}}^{n}\setminus{\mathbb{B}}_{n}(0,1+2\sup_{y\in K}|y|)
\,.
\end{eqnarray}
By inequality (\ref{prop:furasoin2})  and by the elementary inequality
\begin{equation}\label{prop:dora4a}
(|x|-|x-y|)\leq |y|\qquad\forall y\in K\,,
\end{equation}
we have 
\begin{eqnarray}
\label{prop:dora5}
\lefteqn{ 
 \sup_{|x|\geq 1+2\sup_{y\in K}|y|}\,\sup_{y\in K} \biggl\{
e^{ ({\mathrm{Im}}\,k) |x| }|x|^{1+\frac{n-1}{2}}
|\eta_n'(|x-y|)-ik\eta_n(|x-y|)|\biggr\}
} 
\\ \nonumber
&&\qquad\qquad\leq
 \sup_{|x|\geq 1+2\sup_{y\in K}|y|}\,\sup_{y\in K} \biggl\{
e^{({\mathrm{Im}}\,k)(|x|-|x-y|)}
\left(
\frac{|x|}{|x-y|}\right)^{1+\frac{n-1}{2}}
\\ \nonumber
&&\qquad\qquad\quad\times
e^{ ({\mathrm{Im}}\,k) |x-y| }|x-y|^{1+\frac{n-1}{2}}
|\eta_n'(|x-y|)-ik\eta_n(|x-y|)|\biggr\}
\\ \nonumber
&&\qquad \qquad\leq
\sup_{y\in K}e^{({\mathrm{Im}}\,k)|y|}
\left(1+\sup_{y\in K}|y|\right)^{1+\frac{n-1}{2}}
a_2<+\infty\,.
\end{eqnarray}
By inequalities (\ref{prop:furasoin0}), (\ref{prop:dora3})--(\ref{prop:dora4a}) we have 
\begin{eqnarray}
\label{prop:dora6}
\lefteqn{ \sup_{|x|\geq 1+2\sup_{y\in K}|y|}\,\sup_{y\in K} 
\biggl\{e^{ ({\mathrm{Im}}\,k) |x| }|x|^{1+\frac{n-1}{2}}
|k||\eta_n(|x-y|)|\left|\frac{x}{|x|}-\frac{x-y}{|x-y|}\right|\biggr\}
}
\\ \nonumber
&&\qquad \qquad\quad\leq
 \sup_{|x|\geq 1+2\sup_{y\in K}|y|}\,\sup_{y\in K} 
\biggl\{e^{({\mathrm{Im}}\,k)(|x|-|x-y|)}
\left(
\frac{|x|}{|x-y|}\right)^{\frac{n-1}{2}}
\\ \nonumber
&&\qquad \qquad\quad\ \ \times
|k|
\left(
e^{({\mathrm{Im}}\,k) |x-y|}
\eta_n(|x-y|)|x-y|^{\frac{n-1}{2}}
\right)4\sup_{y\in K}|y|\biggr\}
\\ \nonumber
&&\qquad \qquad\leq
\sup_{y\in K}e^{({\mathrm{Im}}\,k)|y|}
\left(1+\sup_{y\in K}|y|\right)^{\frac{n-1}{2}}
|k|a_04\sup_{y\in K}|y|
<+\infty
\,.
\end{eqnarray}
In case $|\zeta|=1$, there exists $r\in \{1,\dots,n\}$ such that $\zeta=e_r$, the $r-th$ unit vector of the canonical basis of ${\mathbb{R}}^n$ and we note that
\begin{eqnarray}
\label{prop:dora7}
\lefteqn{
\frac{x}{|x|} D_{x}\frac{\partial}{\partial x_r}\tilde{S}_{n,k;r}(x-y)-(ik)\frac{\partial}{\partial x_r}\tilde{S}_{n,k;r}(x-y) 
} 
\\ \nonumber
&&
=\sum_{j=1}^n\frac{x_j}{|x|}\biggl\{\frac{\partial}{\partial x_j}\frac{\partial}{\partial x_r}\tilde{S}_{n,k;r}(x-y)\biggr\}-(ik)\frac{\partial}{\partial x_r}\tilde{S}_{n,k;r}(x-y) 
\\ \nonumber
&&
=\sum_{j=1}^n\biggl(\frac{x_j}{|x|}-\frac{x_j-y_j}{|x-y|}\biggr)\frac{\partial}{\partial x_j}\frac{\partial}{\partial x_r}\tilde{S}_{n,k;r}(x-y)
\\ \nonumber
&&
\quad+\sum_{j=1}^n\frac{x_j-y_j}{|x-y|}\frac{\partial}{\partial x_j}\frac{\partial}{\partial x_r}\tilde{S}_{n,k;r}(x-y)
-(ik)\frac{\partial}{\partial x_r}\tilde{S}_{n,k;r}(x-y) \,,
\end{eqnarray}
that
\[
\frac{\partial}{\partial x_r}\tilde{S}_{n,k;r}(x-y)=\eta_n'(|x-y|)\frac{x_r-y_r}{|x-y|}
\]
and that
\begin{eqnarray*}
\lefteqn{
\frac{\partial^2}{\partial x_j\partial x_r}\tilde{S}_{n,k;r}(x-y)=\frac{\partial}{\partial x_j}\frac{\partial}{\partial x_r}\left(\eta_n(x-y)\right)
}
\\ \nonumber
&& 
=\eta_n''(|x-y|)\frac{x_j-y_j}{|x-y|}\frac{x_r-y_r}{|x-y|}+\eta_n'(|x-y|)
\frac{\delta_{rj}|x-y|-(x_r-y_r)\frac{x_j-y_j}{|x-y|}}{|x-y|^2}
\end{eqnarray*}
for all $y\in K$ and $x\in {\mathbb{R}}^{n}\setminus{\mathbb{B}}_{n}(0,1+\sup_{y\in K}|y|)$. Hence, 
  the inequalities (\ref{prop:furasoin1}),  (\ref{prop:furasoin3}) and the elementary inequalities (\ref{prop:dora3})--(\ref{prop:dora4a}) imply that
\begin{eqnarray}
\label{prop:dora8}
\lefteqn{
e^{({\mathrm{Im}}\,k)|x|}|x||x|^{\frac{n-1}{2}}
\biggl|\sum_{j=1}^n\biggl(\frac{x_j}{|x|}-\frac{x_j-y_j}{|x-y|}\biggr)\frac{\partial}{\partial x_j}\frac{\partial}{\partial x_r}\tilde{S}_{n,k;r}(x-y)\biggr|
} 
\\ \nonumber
&&
\leq |x|4 \frac{\sup_{y\in K}|y|}{|x|} \sup_{y\in K} e^{ {\mathrm{Im}}\,k|y|}(1+\sup_{y\in K}|y|)^{ \frac{n-1}{2} }
\biggl\{
a_3+a_1\frac{2}{|x-y|}
\biggr\}
\\ \nonumber
&&
\leq  4  \sup_{y\in K}|y|  \sup_{y\in K} e^{ {\mathrm{Im}}\,k|y|}(1+\sup_{y\in K}|y|)^{ \frac{n-1}{2} }
\biggl\{
a_3+2a_1 
\biggr\}<+\infty
\end{eqnarray}
for all $y\in K$ and $x\in {\mathbb{R}}^{n}\setminus{\mathbb{B}}_{n}(0,1+2\sup_{y\in K}|y|)$. Moreover, 
\begin{eqnarray*}
\lefteqn{
\sum_{j=1}^n\frac{x_j-y_j}{|x-y|}\frac{\partial}{\partial x_j}\frac{\partial}{\partial x_r}\tilde{S}_{n,k;r}(x-y)
-(ik)\frac{\partial}{\partial x_r}\tilde{S}_{n,k;r}(x-y)
}
\\ \nonumber
&&\qquad
=\biggl[\eta_n''(|x-y|)-ik\eta_n'(|x-y|)\biggr] \frac{x_r-y_r}{|x-y|}
\\ \nonumber
&&\qquad\quad
+\eta_n'(|x-y|)\sum_{j=1}^n\frac{x_j-y_j}{|x-y|}
\biggl(\frac{\delta_{rj}}{|x-y|}-\frac{(x_r-y_r)(x_j-y_j)}{|x-y|^3}\biggr)
\end{eqnarray*}
and accordingly, inequalities (\ref{prop:furasoin1}),  (\ref{prop:furasoin4}) and the elementary inequalities (\ref{prop:dora3})--(\ref{prop:dora4a}) imply that
\begin{eqnarray}
\label{prop:dora8a}
\lefteqn{
e^{ ({\mathrm{Im}}\,k) |x| }|x|^{1+\frac{n-1}{2}}\biggl|
\sum_{j=1}^n\frac{x_j-y_j}{|x-y|}\frac{\partial}{\partial x_j}\frac{\partial}{\partial x_r}\tilde{S}_{n,k;r}(x-y)
} 
\\ \nonumber
&& \quad
-(ik)\frac{\partial}{\partial x_r}\tilde{S}_{n,k;r}(x-y)\biggr|
\\ \nonumber
&& 
\leq \sup_{y\in K} e^{{\mathrm{Im}}\,k|y|}(1+\sup_{y\in K}|y|)^{\frac{n-1}{2}+1}
\biggl\{ a_4+2a_1 \biggr\}
\end{eqnarray}
for all $y\in K$ and $x\in {\mathbb{R}}^{n}\setminus{\mathbb{B}}_{n}(0,1+2\sup_{y\in K}|y|)$.  In case $|\zeta|=2$, there exist  $r,s\in \{1,\dots,n\}$ such that $\zeta=e_r+e_s$,   and we note that
\begin{eqnarray}
\label{prop:dora9}
\lefteqn{
\frac{x}{|x|} D_{x}\frac{\partial^2}{\partial x_s\partial x_r}\tilde{S}_{n,k;r}(x-y)-(ik)\frac{\partial^2}{\partial x_s\partial x_r}\tilde{S}_{n,k;r}(x-y) 
} 
\\ \nonumber
&&
=\sum_{j=1}^n\frac{x_j}{|x|}\biggl\{\frac{\partial}{\partial x_j}\frac{\partial^2}{\partial x_s\partial x_r}\tilde{S}_{n,k;r}(x-y)\biggr\}-(ik)\frac{\partial^2}{\partial x_s\partial x_r}\tilde{S}_{n,k;r}(x-y) 
\\ \nonumber
&&
=\sum_{j=1}^n\biggl(\frac{x_j}{|x|}-\frac{x_j-y_j}{|x-y|}\biggr)\frac{\partial}{\partial x_j}\frac{\partial^2}{\partial x_s\partial x_r}\tilde{S}_{n,k;r}(x-y)
\\ \nonumber
&&
\quad+\sum_{j=1}^n\frac{x_j-y_j}{|x-y|}\frac{\partial}{\partial x_j}\frac{\partial^2}{\partial x_s\partial x_r}\tilde{S}_{n,k;r}(x-y)
-(ik)\frac{\partial^2}{\partial x_s\partial x_r}\tilde{S}_{n,k;r}(x-y) \,,
\end{eqnarray}
that
\begin{eqnarray*}
\lefteqn{
\frac{\partial^2}{\partial x_s\partial x_r}\tilde{S}_{n,k;r}(x-y)=\frac{\partial}{\partial x_s}\frac{\partial}{\partial x_r}\left(\eta_n(x-y)\right)
}
\\ \nonumber
&& 
=\eta_n''(|x-y|)\frac{x_s-y_s}{|x-y|}\frac{x_r-y_r}{|x-y|}+\eta_n'(|x-y|)
\frac{\delta_{rs}|x-y|-(x_r-y_r)\frac{x_s-y_s}{|x-y|}}{|x-y|^2}
\end{eqnarray*}
and that
\begin{eqnarray*}
\lefteqn{
\frac{\partial}{\partial x_j}\frac{\partial^2}{\partial x_s\partial x_r}\tilde{S}_{n,k;r}(x-y)
}
\\ \nonumber
&& 
=\eta_n'''(|x-y|)\frac{x_j-y_j}{|x-y|}\frac{x_s-y_s}{|x-y|}\frac{x_r-y_r}{|x-y|}
\\ \nonumber
&&\quad
+\eta_n''(|x-y|)\biggl(\frac{\delta_{sj}}{|x-y|}-\frac{(x_s-y_s)(x_j-y_j)}{|x-y|^3}\biggr)\frac{x_r-y_r}{|x-y|}
 \\ \nonumber
&&\quad
+\eta_n''(|x-y|)\frac{x_s-y_s}{|x-y|}\biggl(\frac{\delta_{rj}}{|x-y|}-\frac{(x_r-y_r)(x_j-y_j)}{|x-y|^3}\biggr) 
 \\ \nonumber
&&\quad
+\eta_n''(|x-y|)\frac{x_j-y_j}{|x-y|}\biggl(\frac{\delta_{rs}}{|x-y|}-\frac{(x_r-y_r)(x_s-y_s)}{|x-y|^3}\biggr) 
\\ \nonumber
&&\quad
+\eta_n'(|x-y|)
\biggl(
-\delta_{rs}\frac{x_j-y_j}{|x-y|^3}
-\frac{[\delta_{rj}(x_s-y_s)+(x_r-y_r)\delta_{sj}]}{|x-y|^3}
\\ \nonumber
&&\quad
+\frac{(x_r-y_r)(x_s-y_s)3(x_j-y_j) 
}{|x-y|^5}
\biggr)
\end{eqnarray*}
for all $y\in K$ and $x\in {\mathbb{R}}^{n}\setminus{\mathbb{B}}_{n}(0,1+\sup_{y\in K}|y|)$. Hence, 
  the inequalities (\ref{prop:furasoin1}),  (\ref{prop:furasoin3}),  (\ref{prop:furasoin5}) and the elementary inequalities (\ref{prop:dora3})--(\ref{prop:dora4a}) imply that
\begin{eqnarray}
\label{prop:dora10}
\lefteqn{e^{({\mathrm{Im}}\,k)|x|}\,|x||x|^{\frac{n-1}{2}}
\biggl|\sum_{j=1}^n\biggl(\frac{x_j}{|x|}-\frac{x_j-y_j}{|x-y|}\biggr)\frac{\partial}{\partial x_j}\frac{\partial^2}{\partial x_s\partial x_r}\tilde{S}_{n,k;r}(x-y)\biggr|
} 
\\ \nonumber
&&\leq
|x|4 \frac{\sup_{y\in K}|y|}{|x|} \sup_{y\in K} e^{ {\mathrm{Im}}\,k|y|}(1+\sup_{y\in K}|y|)^{ \frac{n-1}{2} }
\biggl\{
a_5+\frac{6 a_3}{|x-y|}+\frac{6a_1}{|x-y|^2}
\biggr\}
\\ \nonumber
&&\leq
 4  \sup_{y\in K}|y| \sup_{y\in K} e^{ {\mathrm{Im}}\,k|y|}(1+\sup_{y\in K}|y|)^{ \frac{n-1}{2} }
\biggl\{
a_5+ 6 a_3 + 6a_1 
\biggr\}<+\infty
\end{eqnarray}
for all $y\in K$ and $x\in {\mathbb{R}}^{n}\setminus{\mathbb{B}}_{n}(0,1+2\sup_{y\in K}|y|)$. Moreover, 
\begin{eqnarray*}
\lefteqn{
\sum_{j=1}^n\frac{x_j-y_j}{|x-y|}\frac{\partial}{\partial x_j}\frac{\partial^2}{\partial x_s\partial x_r}\tilde{S}_{n,k;r}(x-y)
-(ik)\frac{\partial^2}{\partial x_s\partial x_r}\tilde{S}_{n,k;r}(x-y)
}
\\ \nonumber
&&\qquad
=\biggl[\eta_n'''(|x-y|)-ik\eta_n''(|x-y|)\biggr]
 \frac{x_r-y_r}{|x-y|}\frac{x_s-y_s}{|x-y|}
\\ \nonumber
&&\qquad\quad
+\eta_n''(|x-y|)\sum_{j=1}^n\frac{x_j-y_j}{|x-y|}\biggl[\frac{\delta_{sj}}{|x-y|}-\frac{(x_s-y_s)(x_j-y_j)}{|x-y|^3}\biggr]\frac{x_r-y_r}{|x-y|}
\\ \nonumber
&&\qquad\quad
+\eta_n''(|x-y|)\frac{x_s-y_s}{|x-y|}\sum_{j=1}^n\frac{x_j-y_j}{|x-y|}
\biggl[\frac{\delta_{rj}}{|x-y|}-\frac{(x_r-y_r)(x_j-y_j)}{|x-y|^3}\biggr]
\\ \nonumber
&&\qquad\quad
+\eta_n''(|x-y|)\sum_{j=1}^n\frac{x_j-y_j}{|x-y|}\frac{x_j-y_j}{|x-y|}
\biggl[\frac{\delta_{rs}}{|x-y|}-\frac{(x_r-y_r)(x_s-y_s)}{|x-y|^3}\biggr]
\\ \nonumber
&&\qquad\quad
+\eta_n'(|x-y|)\biggl[
-\delta_{rs}\sum_{j=1}^n\frac{x_j-y_j}{|x-y|}\frac{x_j-y_j}{|x-y|^3}
\\ \nonumber
&&\qquad\quad
-\sum_{j=1}^n\frac{x_j-y_j}{|x-y|}\frac{\delta_{rj}(x_s-y_s)+(x_r-y_r)\delta_{sj}}{|x-y|^3}
\\ \nonumber
&&\qquad\quad
+\sum_{j=1}^n\frac{x_j-y_j}{|x-y|}\frac{(x_r-y_r)(x_s-y_s)3 (x_j-y_j)}{|x-y|^5}
\biggl]
\\ \nonumber
&&\qquad\quad
-ik\eta_n'(|x-y|)\biggl[
\frac{\delta_{rs}}{|x-y|}-\frac{(x_r-y_r)(x_s-y_s)}{|x-y|^3}
\biggr]
\end{eqnarray*}
and accordingly  the inequalities (\ref{prop:furasoin1}),  (\ref{prop:furasoin3}),  (\ref{prop:furasoin6}) and the elementary inequalities (\ref{prop:dora3})--(\ref{prop:dora4a}) imply that
\begin{eqnarray}
\label{prop:dora10a}
\lefteqn{
e^{ ({\mathrm{Im}}\,k) |x| }|x|^{1+\frac{n-1}{2}}\biggl|
\sum_{j=1}^n\frac{x_j-y_j}{|x-y|}\frac{\partial}{\partial x_j}\frac{\partial^2}{\partial x_s\partial x_r}\tilde{S}_{n,k;r}(x-y)
} 
\\ \nonumber
&& \quad
-(ik)\frac{\partial^2}{\partial x_s\partial x_r}\tilde{S}_{n,k;r}(x-y)\biggr|
\\ \nonumber
&& 
\leq \sup_{y\in K} e^{{\mathrm{Im}}\,k|y|}(1+\sup_{y\in K}|y|)^{\frac{n-1}{2}}
\biggl\{
\frac{|x|}{|x-y|}a_6+|x|a_3\frac{2}{|x-y|}+|x|a_3\frac{2}{|x-y|}
\\ \nonumber
&& \quad
+|x|a_3\frac{2}{|x-y|}
+|x|a_1\biggl[\frac{1}{|x-y|^2} +\frac{2}{|x-y|^2}+\frac{3}{|x-y|^2}\biggr]
\\ \nonumber
&& \quad
+
|k|a_1|x|\biggl[
\frac{1}{|x-y|}+\frac{1}{|x-y|}
\biggr]
\biggr\}
\\ \nonumber
&& 
\leq \sup_{y\in K} e^{{\mathrm{Im}}\,k|y|}(1+\sup_{y\in K}|y|)^{\frac{n+1}{2}}
\biggl\{a_6+6a_3+a_1\frac{6}{|x-y|}+2|k|a_1\biggr\}<+\infty
\end{eqnarray}
for all $y\in K$ and $x\in {\mathbb{R}}^{n}\setminus{\mathbb{B}}_{n}(0,1+2\sup_{y\in K}|y|)$. 
 
By combining inequalities (\ref{prop:dora2}), (\ref{prop:dora5}), (\ref{prop:dora6}), equality (\ref{prop:dora7}), inequalities (\ref{prop:dora8}), (\ref{prop:dora8a}), equality (\ref{prop:dora9}), inequalities (\ref{prop:dora10}), (\ref{prop:dora10a}), we conclude that inequality (\ref{prop:dora0}) holds true. We now prove inequality   (\ref{prop:dora1}). First we note that
\[
\frac{\partial}{\partial v(y)}\tilde{S}_{n,k;r}(x-y)=
\frac{\partial}{\partial v(y)}\left(\eta_n (|x-y|)\right)=
-\eta_n'(|x-y|)\frac{x-y}{|x-y|}\cdot v(y)\,,
\]
for all $x\in {\mathbb{R}}^{n}\setminus K$ and $y\in K$, and that
\begin{eqnarray}
\label{prop:dora12}
\lefteqn{
\frac{x}{|x|}\cdot D_{x}\left\{
\frac{\partial}{\partial v(y)}\eta_n (|x-y|)
\right\}
-ik
\frac{\partial}{\partial v(y)}\eta_n (|x-y|)
}
\\ \nonumber
&& =
-\frac{x}{|x|}\cdot D_{x}\left\{
\eta_n'(|x-y|)\frac{x-y}{|x-y|}\cdot v(y)
\right\}
+ik\eta_n'(|x-y|)\frac{x-y}{|x-y|}\cdot v(y)
\\ \nonumber
&& =
-\frac{x}{|x|}\cdot\biggl\{\biggr.
\eta_n''(|x-y|)\frac{x-y}{|x-y|}\biggl(
\frac{x-y}{|x-y|}\cdot v(y)
\biggr)
\\ \nonumber
&&\quad
+\eta_n'(|x-y|)\biggl[
\frac{v(y)}{|x-y|}
-
(x-y)\frac{(x-y)\cdot v(y)}{|x-y|^{3}}
\biggr]\biggl.\biggr\}
\\ \nonumber
&&\quad
+ik\eta_n'(|x-y|)\frac{(x-y)\cdot v(y)}{|x-y|}
\\ \nonumber
&& =
-\frac{x}{|x|}\cdot\biggl\{
\eta_n''(|x-y|)\biggl(\frac{x-y}{|x-y|}
-\frac{x}{|x|}\biggr)\biggl(
\frac{x-y}{|x-y|}\cdot v(y)
\biggr)\biggr\}
\\ \nonumber
&&\quad
-\frac{x}{|x|}\cdot\biggl\{
\eta_n''(|x-y|)\biggl( \frac{x}{|x|}\biggr)\biggl(
\frac{x-y}{|x-y|}\cdot v(y)
\biggr)\biggr\}
\\ \nonumber
&&\quad
-\eta_n'(|x-y|)\frac{x}{|x|} \cdot
\biggl[
\frac{v(y)}{|x-y|}
-
(x-y)\frac{(x-y)\cdot v(y)}{|x-y|^{3}}
\biggr]
\\ \nonumber
&&\quad
+ik\eta_n'(|x-y|)\left(
\frac{(x-y)\cdot v(y)}{|x-y|}
\right)
\\ \nonumber
&& =
-\frac{x}{|x|}\cdot\biggl\{
\eta_n''(|x-y|)\biggl(\frac{x-y}{|x-y|}
-\frac{x}{|x|}\biggr)\biggl(
\frac{x-y}{|x-y|}\cdot v(y)
\biggr)\biggr\}
\\ \nonumber
&&\quad
- \eta_n'(|x-y|)\frac{x}{|x|}\cdot
\biggl[
\frac{v(y)}{|x-y|}
-
(x-y)\frac{(x-y)\cdot v(y)}{|x-y|^{3}}
\biggr]
\\ \nonumber
&&\quad
-\biggl\{
\eta_n''(|x-y|)-ik\eta_n'(|x-y|)
\biggr\}
\frac{(x-y)\cdot v(y)}{|x-y| }
\,,
\end{eqnarray}
for all $y\in K$ and   $x\in {\mathbb{R}}^{n}\setminus{\mathbb{B}}_{n}(0,1+ \sup_{y\in K}|y|)$.
By inequality (\ref{prop:furasoin3})  and inequalities (\ref{prop:dora3})--(\ref{prop:dora4a}), we obtain
\begin{eqnarray}
\label{prop:dora13}
\lefteqn{
 \sup_{  |x|\geq 1+2\sup_{y\in K}|y|   }\,
 \sup_{y\in K} 
  e^{ ({\mathrm{Im}}\,k)|x|  }|x||x|^{\frac{n-1}{2}}
 }
 \\ \nonumber
&&\quad\times
 \biggl|
 \frac{x}{|x|}\cdot
 \biggl\{
\eta_n''(|x-y|)\biggl(
\frac{x-y}{|x-y|}
-\frac{x}{|x|}
\biggr)
\biggl(
\frac{x-y}{|x-y|}\cdot v(y)
\biggr)
\biggr\}
\biggr|
\\ \nonumber
&&
\leq
\sup_{|x|\geq 1+2\sup_{y\in K}|y|}\,\sup_{y\in K} 
e^{ ({\mathrm{Im}}\,k)[|x|-|x-y|] }
\left(
\frac{|x|}{|x-y|}
\right)^{ \frac{n-1}{2}  }
e^{  ({\mathrm{Im}}\,k)|x-y|  }  |x-y|^{ \frac{n-1}{2}  }
\\ \nonumber
&&\quad\times
|\eta_n''(|x-y|)|4\sup_{y\in K} |y|<+\infty\,.
\end{eqnarray}
Then by inequality (\ref{prop:furasoin1})  and inequalities (\ref{prop:dora3}), (\ref{prop:dora4a}), we obtain
\begin{eqnarray}
\label{prop:dora14}
\lefteqn{
 \sup_{  |x|\geq 1+2\sup_{y\in K}|y|   }\,
 \sup_{y\in K} 
  e^{ ({\mathrm{Im}}\,k)|x|  }|x||x|^{\frac{n-1}{2}}
 }
 \\ \nonumber
&&\quad\times
\biggl|
-\frac{x}{|x|}\eta_n'(|x-y|)
\biggl[
\frac{v(y)}{|x-y|}
-
(x-y)\frac{(x-y)\cdot v(y)}{|x-y|^{3}}
\biggr]
\biggr|
\\ \nonumber
&&
\leq
\sup_{|x|\geq 1+2\sup_{y\in K}|y|}\,\sup_{y\in K} 
e^{ ({\mathrm{Im}}\,k)[|x|-|x-y|] }
\left(
\frac{|x|}{|x-y|}
\right)^{ \frac{n-1}{2}  }
\\ \nonumber
&&\quad\times
e^{  ({\mathrm{Im}}\,k) |x-y|  }  |x-y|^{ \frac{n-1}{2}  }
|\eta_n'(|x-y|)|\frac{2|x|}{|x-y|}<+\infty\,.
\end{eqnarray}
Then by inequality (\ref{prop:furasoin4}) and inequalities (\ref{prop:dora3}), (\ref{prop:dora4a}), we obtain
\begin{eqnarray}
\label{prop:dora15}
\lefteqn{
 \sup_{  |x|\geq 1+2\sup_{y\in K}|y|   }\,
 \sup_{y\in K} 
  e^{ ({\mathrm{Im}}\,k)|x|  }|x||x|^{\frac{n-1}{2}}
 }
 \\ \nonumber
&&\quad\times
\left|
\biggl\{
\eta_n''(|x-y|)-ik\eta_n'(|x-y|)
\biggr\}
\frac{(x-y)\cdot v(y)}{|x-y| }
\right|
\\ \nonumber
&&
\leq
\sup_{|x|\geq 1+2\sup_{y\in K}|y|}\,\sup_{y\in K} 
e^{ ({\mathrm{Im}}\,k)[|x|-|x-y|] }
\left(
\frac{|x|}{|x-y|}
\right)^{ \frac{n-1}{2} +1 }
\\ \nonumber
&&\quad\times
|x-y|^{ \frac{n-1}{2} +1 }
e^{ ({\mathrm{Im}}\,k)|x-y|}
|\eta_n''(|x-y|)-ik\eta_n'(|x-y|)|<+\infty\,.
\end{eqnarray}
By equality (\ref{prop:dora12}) and by inequalities (\ref{prop:dora13}), (\ref{prop:dora14}), (\ref{prop:dora15}), we deduce the validity of inequality (\ref{prop:dora1}).
\hfill  $\Box$ 

\vspace{\baselineskip}

 Now we check that the layer potentials corresponding to the fundamental solution that satisfies the outgoing $(k)$-radiation condition also satisfy 
the outgoing $(k)$-radiation condition. For the double layer potential the statement is classic and well-known (cf.~\textit{e.g.}, Colton and Kress \cite[Thm.~3.2,]{CoKr92}). Instead,  for the single layer potential
with densities in the space $V^{-1,\alpha}(\partial\Omega)$, the author could find no reference in the literature.
\begin{theorem}
\label{thm:pora}
Let  $k\in {\mathbb{C}}\setminus ]-\infty,0]$, ${\mathrm{Im}}\,k\geq 0$, $\alpha\in]0,1[$. Let $\Omega$ be a bounded open subset of ${\mathbb{R}}^{n}$ of class $C^{1,\alpha}$. Then the following statements hold.
 \begin{enumerate}
\item[(i)] Let $\mu\in V^{-1,\alpha}(\partial\Omega)$. Then $v_\Omega[\tilde{S}_{n,k;r},\mu]$ satisfies the outgoing $(k)$-radiation condition. More precisely,
\begin{equation}
\label{thm:pora2}
\limsup_{x\to\infty}|x|e^{{\mathrm{Im}}\, k |x|}|x|^{\frac{n-1}{2}}
\left(
\frac{x}{|x|}Dv_\Omega[\tilde{S}_{n,k;r},\mu](x)-ik
v_\Omega[\tilde{S}_{n,k;r},\mu](x)
\right)<+\infty\,.
\end{equation}
\item[(ii)] Let $\mu\in L^{1}(\partial\Omega)$. Then  $w_\Omega[\tilde{S}_{n,k;r},\mu]$ 
satisfies the outgoing $(k)$-radiation condition. More precisely,
\begin{equation}
\label{thm:pora4}
\limsup_{x\to\infty}|x|e^{{\mathrm{Im}}\, k |x|}|x|^{\frac{n-1}{2}}
\left(
\frac{x}{|x|}Dw_\Omega[\tilde{S}_{n,k;r},\mu](x)-ik
w_\Omega[\tilde{S}_{n,k;r},\mu](x)
\right)<+\infty\,.
\end{equation}
\end{enumerate}
\end{theorem}
{\bf Proof.} (i) Since $\tilde{S}_{n,k;r}\in C^\infty({\mathbb{R}}^n\setminus\{0\})$, and 
$\tilde{S}_{n,k;r}$, $\frac{\partial \tilde{S}_{n,k;r}}{\partial x_j}$ for $j\in\{1,\dots,n\}$ are locally integrable in ${\mathbb{R}}^n$, classical properties of the convolution imply that
\[
\frac{\partial v_\Omega[\tilde{S}_{n,k;r},\mu]}{\partial x_j}
=(r_{|\partial\Omega}^t\mu)\ast \frac{\partial \tilde{S}_{n,k;r}}{\partial x_j}
=\langle r_{|\partial\Omega}^t\mu(y), \frac{\partial \tilde{S}_{n,k;r}}{\partial x_j}(\cdot-y)\rangle 
\qquad\text{in}\ {\mathbb{R}}^n\setminus\partial\Omega\,,
\]
 for all $j\in\{1,\dots,n\}$.  Then Proposition \ref{prop:dora} implies that
\begin{eqnarray}
\label{thm:pora5}
\lefteqn{
\left|
|x|^{\frac{n-1}{2}}
\left(
\frac{x}{|x|}Dv_\Omega[\tilde{S}_{n,k;r},\mu](x)-ik
v_\Omega[\tilde{S}_{n,k;r},\mu](x)
\right)
\right|
}
\\ \nonumber
&&
=|x|^{\frac{n-1}{2}}
\biggl|
\sum_{j=1}^n\frac{x_j}{|x|}\frac{\partial}{\partial x_j}
\langle\mu(y),\tilde{S}_{n,k;r}(x-y)\rangle
\\ \nonumber
&&\quad
-ik 
\langle\mu(y),\tilde{S}_{n,k;r} (x-y) \rangle
\biggr|
\\ \nonumber
&&
=|x|^{\frac{n-1}{2}}
\biggl|
 \langle\mu(y),\sum_{j=1}^n\frac{x_j}{|x|}\frac{\partial}{\partial x_j}\tilde{S}_{n,k;r}(x-y)\rangle
 \\ \nonumber
&&\quad
-ik 
\langle\mu(y),\tilde{S}_{n,k;r} (x-y) \rangle
\biggr|
\\ \nonumber
&&
=
 |x|^{\frac{n-1}{2}}
\biggl|\langle\mu(y),
\frac{x}{|x|}D\tilde{S}_{n,k;r}(x-y)-ik\tilde{S}_{n,k;r}(x-y)
\rangle\biggr|
\\ \nonumber
&&
\leq \|\mu\|_{V^{-1,\alpha}(\partial\Omega)}|x|^{\frac{n-1}{2}}
\\ \nonumber
&&\quad\times
\left\|
\frac{x}{|x|}D\tilde{S}_{n,k;r}(x-\cdot)-ik\tilde{S}_{n,k;r}(x-\cdot)
\right\|_{C^{1,\alpha}(\partial\Omega)}
\\ \nonumber
&&
\leq \|\mu\|_{V^{-1,\alpha}(\partial\Omega)}|x|^{\frac{n-1}{2}}\iota
\\ \nonumber
&&\quad\times
\left\|
\frac{x}{|x|}D\tilde{S}_{n,k;r}(x-\cdot)-ik\tilde{S}_{n,k;r}(x-\cdot)
\right\|_{C^{2}(\overline{{\mathbb{B}}_n(0,\sup_{y\in \partial\Omega}|y|)})}
\\ \nonumber
&&
\leq
c_{0,\overline{{\mathbb{B}}_n(0,\sup_{y\in \partial\Omega}|y|)}}\|\mu\|_{V^{-1,\alpha}(\partial\Omega)}
|x|^{-1}e^{-{\mathrm{Im}}\,k|x|}\,,
 \end{eqnarray}
for all $x\in {\mathbb{R}}^{n}\setminus {\mathbb{B}}_{n} (0,1+2
\sup_{y\in \partial\Omega}|y|)$. Here $\iota$ denotes the norm of the restriction map from $C^{2}(\overline{{\mathbb{B}}_n(0,\sup_{y\in \partial\Omega}|y|)})$ to $C^{1,\alpha}(\partial\Omega)$. Then inequality (\ref{thm:pora5}) implies the validity of statement  (i). 

We now consider statement (ii). By Proposition \ref{prop:dora} and by a classical result on the differentiation of integrals depending on a parameter, we have 
\begin{eqnarray}
\label{thm:pora6}
\lefteqn{
\left|
|x|^{\frac{n-1}{2}}
\left(
\frac{x}{|x|}Dw_\Omega[\tilde{S}_{n,k;r},\mu](x)-ik
w_\Omega[\tilde{S}_{n,k;r},\mu](x)
\right)
\right|
}
\\ \nonumber
&&
\leq
\biggl|\biggr.
|x|^{\frac{n-1}{2}}
\int_{\partial\Omega}
\biggl[\biggr.
\frac{x}{|x|}D_{x}\left(
\frac{\partial}{\partial\nu_{\Omega}(y)}
\tilde{S}_{n,k;r}(x-y)\right)
\\ \nonumber
&&\qquad\qquad\qquad\qquad
-ik\frac{\partial}{\partial\nu_{\Omega}(y)}\tilde{S}_{n,k;r}(x-y)
\biggl.\biggr]\mu(y)\,d\sigma_{y}
\biggl.\biggr|
\\ \nonumber
&&
\leq
c_{1,\partial\Omega}\int_{\partial\Omega}|\mu|\,d\sigma_{y}
 |x|^{-1}e^{-{\mathrm{Im}}\,k|x|}\,,
 \end{eqnarray}
for all $x\in {\mathbb{R}}^{n}\setminus {\mathbb{B}}_{n} (0,1+2
\sup_{y\in \partial\Omega}|y|)$. Then 
Inequality (\ref{thm:pora6}) implies the validity of statement (ii). 
\hfill  $\Box$ 

\vspace{\baselineskip}

\section{Preliminaries to the acoustic volume potential}
 \label{sec:acvolpot}
Let $\alpha\in]0,1]$ and $m\in {\mathbb{N}}$.  If $\Omega$ is a  bounded open subset of ${\mathbb{R}}^{n}$, then we can consider the restriction map $r_{|\overline{\Omega}}$ from ${\mathcal{D}}({\mathbb{R}}^n)$ to $C^{m,\alpha}(\overline{\Omega})$. Then the transpose map $r_{|\overline{\Omega}}^t$ is linear and continuous from $(C^{m,\alpha}(\overline{\Omega}))'$ to ${\mathcal{D}}'({\mathbb{R}}^n)$. Moreover, if $\mu\in (C^{m,\alpha}(\overline{\Omega}))'$, then $r_{|\overline{\Omega}}^t\mu$ has compact support. Hence, it makes sense to consider the convolution of 
  $r_{|\overline{\Omega}}^t\mu$ with  the fundamental solution of either the Laplace or the Helmholtz  operator. Thus we are now ready to introduce the following known definition.
 \begin{definition}\label{defn:dvpsl}
 Let $\alpha\in]0,1]$, $m\in {\mathbb{N}}$. 
 Let   $\Omega$ be a bounded open subset of ${\mathbb{R}}^{n}$. Let $\lambda\in {\mathbb{C}}$.  Let $S_{n,\lambda} $ be a fundamental solution 
of the operator $\Delta+\lambda$. If $\mu\in (C^{m,\alpha}(\overline{\Omega}))'$, then the (distributional) volume potential relative to $S_{n,\lambda} $ and $\mu$ is the distribution
\[
{\mathcal{P}}_\Omega[S_{n,\lambda} ,\mu]=(r_{|\overline{\Omega}}^t\mu)\ast S_{n,\lambda}  \in {\mathcal{D}}'({\mathbb{R}}^n)\,.
\]
\end{definition}
Under the assumptions of Definition \ref{defn:dvpsl}, we  set
\begin{eqnarray}\label{prop:dvpsl3}
{\mathcal{P}}_\Omega^+[S_{n,\lambda} ,\mu]&\equiv&\left((r_{|\overline{\Omega}}^t\mu)\ast S_{n,\lambda}  \right)_{|\Omega}
\qquad\text{in}\ \Omega\,,
\\ \nonumber
{\mathcal{P}}_\Omega^-[S_{n,\lambda} ,\mu]  &\equiv&
\left((r_{|\overline{\Omega}}^t\mu)\ast S_{n,\lambda}  \right)_{|\Omega^-}
\qquad\text{in}\ \Omega^-\,.
\end{eqnarray}
In general, $(r_{|\overline{\Omega}}^t\mu)\ast S_{n,\lambda} $ is not a function, \textit{i.e.} $(r_{|\overline{\Omega}}^t\mu)\ast S_{n,\lambda} $ is not a distribution that is associated to a locally integrable function in ${\mathbb{R}}^n$. 
However, this is the case if for example  $\mu$ is associated to a function $f$ of $ L^\infty(\Omega)$,  and thus 
 the (distributional) volume potential relative to $S_{n,\lambda} $ and $\mu$ is associated to the function
\begin{equation}\label{prop:dvpsa1}
\int_{\Omega}S_{n,\lambda} (x-y)f(y)\,dy\qquad{\mathrm{a.a.}}\ x\in {\mathbb{R}}^n\,,
\end{equation}
that is locally integrable in ${\mathbb{R}}^n$  and that with some abuse of notation we  denote by the symbol    ${\mathcal{P}}_\Omega[S_{n,\lambda} , f]$.  Then the following classical result is known (cf.~\textit{e.g.}, \cite[\S 6]{La24d}).
\begin{theorem}\label{thm:nwtdma} 
 Let $m\in {\mathbb{N}}$, $\alpha\in]0,1[$. Let $\Omega$ be a bounded open subset of  ${\mathbb{R}}^n$ of class $C^{m+1,\alpha}$. Let $\lambda\in {\mathbb{C}}$.  Let $S_{n,\lambda}  $ be a fundamental solution 
of the operator $\Delta+\lambda$. Then the following statements hold. 
 
 \item[(i)]  ${\mathcal{P}}_\Omega^+[S_{n,\lambda} ,\cdot]$ is linear and continuous from $C^{m,\alpha}(\overline{\Omega})$ to $C^{m+2,\alpha}(\overline{\Omega})$.
\item[(ii)]   ${\mathcal{P}}_\Omega^-[S_{n,\lambda} ,\cdot]$ is linear and continuous from $C^{m,\alpha}(\overline{\Omega})$ to   $C^{m+2,\alpha}(\overline{{\mathbb{B}}_n(0,r)}\setminus\Omega)$ for all $r\in]0,+\infty[$ such that $\overline{\Omega}\subseteq {\mathbb{B}}_n(0,r)$.
\end{theorem}
 Instead, for a Schauder space of negative exponent, the following statement holds (cf.~\textit{e.g.}, \cite[\S 7]{La24d}), that is a  generalization to volume potentials of  nonhomogeneous second order elliptic operators of a known result for the Laplace operator  (cf.~\cite[Thm.~3.6 (ii)]{La08a}, Dalla Riva, the author and Musolino~\cite[Thm.~7.19]{DaLaMu21}).
    \begin{proposition}\label{prop:dvpsnecr-1a}
  Let $\alpha\in]0,1[$.    Let $\Omega$ be a bounded open subset of ${\mathbb{R}}^{n}$ of class $C^{1,\alpha}$. Let $\lambda\in {\mathbb{C}}$.  Let $S_{n,\lambda}  $ be a fundamental solution 
of the operator $\Delta+\lambda$.   Then the following statements hold. 
 \begin{enumerate}
\item[(i)] If  $f=  f_{0}+\sum_{j=1}^{n}\frac{\partial}{\partial x_{j}}f_{j}\in C^{-1,\alpha}(\overline{\Omega}) $, then
 \begin{equation}\label{prop:dvpsnecr-1a2}
{\mathcal{P}}_\Omega^+[S_{n,\lambda} ,E^\sharp[f]]\in C^{1,\alpha}(\overline{\Omega}), \ 
{\mathcal{P}}_\Omega^-[S_{n,\lambda} ,E^\sharp[f]]\in C^{1,\alpha}_{{\mathrm{loc}} }(\overline{\Omega^-}) \,,\end{equation}
and 
\begin{equation}\label{defn:Ppm1}
{\mathcal{P}}_\Omega^+[S_{n,\lambda} ,E^\sharp[f]](x)={\mathcal{P}}_\Omega^-[S_{n,\lambda} ,E^\sharp[f]](x)\qquad\forall x\in\partial\Omega\,.
\end{equation}
Moreover,
\begin{eqnarray}\label{prop:dvpsnecr-1a2a}
&&\Delta {\mathcal{P}}_\Omega^+[S_{n,\lambda} ,E^\sharp[f]] 
+\lambda{\mathcal{P}}_\Omega^+[S_{n,\lambda} ,E^\sharp[f]]= f\qquad\textit{in}\ {\mathcal{D}}'(\Omega)\,,
\\ \nonumber
&&\Delta  {\mathcal{P}}_\Omega^-[S_{n,\lambda} ,E^\sharp[f]]
+\lambda  {\mathcal{P}}_\Omega^-[S_{n,\lambda} ,E^\sharp[f]]= 0\qquad\textit{in}\ {\mathcal{D}}'({\mathbb{R}}^n\setminus\overline{\Omega})
\,.
\end{eqnarray}
\item[(ii)] The linear operator  ${\mathcal{P}}_\Omega^+[S_{n,\lambda} ,E^\sharp[\cdot]]$ is  continuous from $C^{-1,\alpha}(\overline{\Omega}) $ to $C^{1,\alpha}(\overline{\Omega})$.
\item[(iii)] Let $r\in ]0,+\infty[$ be such that $\overline{\Omega}\subseteq {\mathbb{B}}_n(0,r)$. Then  ${\mathcal{P}}_\Omega^-[S_{n,\lambda} ,E^\sharp[\cdot]]_{|\overline{{\mathbb{B}}_n(0,r)}\setminus\Omega}$ is linear and continuous from $C^{-1,\alpha}(\overline{\Omega}) $ to   $C^{1,\alpha}(\overline{{\mathbb{B}}_n(0,r)}\setminus\Omega)$.
\end{enumerate}  
(See Proposition \ref{prop:nschext} for the definition of $E^\sharp$).
\end{proposition}
 \section{A distributional form of the second and third Green Identities}\label{sec:dgridh}
  
\begin{theorem}\label{thm:dsecondgreenh}
 Let  $\alpha\in]0,1[$. Let $\Omega$ be a bounded open subset of ${\mathbb{R}}^n$ of class $C^{1,\alpha}$.  If $u\in C^{0,\alpha}(\overline{\Omega})_\Delta$ and $v\in C^{1,\alpha}(\overline{\Omega})_\Delta$, then
\begin{equation}\label{thm:dsecondgreenh1}
\langle E^\sharp[\Delta u],v\rangle -\int_\Omega u\Delta v\,dx=
\langle \frac{\partial u}{\partial\nu},v\rangle -\int_{\partial\Omega}\frac{\partial v}{\partial\nu}u\,d\sigma\,,
\end{equation}
where $\frac{\partial u}{\partial\nu}$ denotes the distributional normal derivative of Definition \ref{defn:conoderdedu} and $E^\sharp$ is the linear extension operator of Proposition \ref{prop:nschext}.
\end{theorem}
{\bf Proof.} By the known  approximation Lemma \ref{lem:apr0ad}, there exists a sequence $\{u_j\}_{j\in {\mathbb{N}}}$ in $C^{1,\alpha}(\overline{\Omega})_\Delta$ as in (\ref{lem:apr1ad1}). Since $v$, $u_j\in C^{1,\alpha}(\overline{\Omega})_\Delta
\subseteq C^2(\Omega)\cap C^{1,\alpha}(\overline{\Omega})$, the following second Green Identity holds
\[
\int_\Omega v\Delta u_j -u_j\Delta v\,dx=\int_{\partial\Omega}\frac{\partial u_j}{\partial\nu}v-
\frac{\partial v}{\partial\nu}u_j\,d\sigma
\]
(cf.~\textit{e.g.},  \cite[Thm.~4.3]{DaLaMu21} and (\ref{eq:c1adc2})) and accordingly   equalities  (\ref{prop:nschext3}), (\ref{lem:conoderdeducl1}) imply that
\begin{equation}\label{thm:dsecondgreenh2}
\langle E^\sharp[\Delta u_j],v\rangle  -\int_\Omega u_j \Delta v\,dx=\langle \frac{\partial u_j}{\partial\nu},v\rangle -
\int_{\partial\Omega}\frac{\partial v}{\partial\nu}u_j\,d\sigma
\end{equation}
for all $j\in {\mathbb{N}}$. Then the membership of $v\in C^{1,\alpha}(\overline{\Omega})\subseteq C^{1,\beta}(\overline{\Omega})$, the continuity of $E^\sharp$ from $C^{-1,\beta}(\overline{\Omega})$ to
$\left(C^{1,\beta}(\overline{\Omega})\right)'$ and the continuity of the distributional normal derivative from 
$C^{0,\beta}(\overline{\Omega})_\Delta$ to $V^{-1,\beta}(\partial\Omega)\subseteq \left(C^{1,\beta}(\partial\Omega)\right)'$ and the limiting relation of (\ref{lem:apr1ad1}) imply that by taking the limit as $j$ tends to infinity in (\ref{thm:dsecondgreenh2}), we obtain the distributional second Green Identity of the statement.\hfill  $\Box$ 

\vspace{\baselineskip}

Next we introduce the following distributional form of the  third Green Identity.

\begin{theorem}\label{thm:dthirdgreenh}
Let  $\alpha\in]0,1[$. Let $\Omega$ be a bounded open subset of ${\mathbb{R}}^n$ of class $C^{1,\alpha}$. Let $\lambda\in {\mathbb{C}}$.  Let $S_{n,\lambda}$ be a fundamental solution 
of the operator $\Delta+\lambda$. If $u\in C^{0,\alpha}(\overline{\Omega})_\Delta $, then the following statements hold, where $\frac{\partial u}{\partial\nu_{\Omega}}$ denotes the distributional normal derivative of of Definition \ref{defn:conoderdedu} and ${\mathcal{P}}_\Omega^+[S_{ \lambda },E^\sharp[\Delta u+\lambda u]]$ denotes the distributional volume potential of $E^\sharp[\Delta u+\lambda u]$ as in Definition \ref{defn:dvpsl} (see also (\ref{prop:dvpsl3})).
\begin{enumerate}
\item[(i)]
\begin{eqnarray}
\label{thm:dthirdgreenh1}
\lefteqn{
u(x)={\mathcal{P}}_\Omega^+[S_{ \lambda },E^\sharp[\Delta u+\lambda u]](x)
}
\\ \nonumber
&&  
+\int_{\partial\Omega}
u(y)\frac{\partial}{\partial\nu_{\Omega,y} }\left(S_{n,\lambda} (x-y)\right)\,d\sigma_{y}
-\langle \frac{\partial u}{\partial\nu_{\Omega}}(y),S_{n,\lambda} (x-y)\rangle 
\quad\forall x\in\Omega\,,
\\ \label{thm:dthirdgreenh1a}
\lefteqn{
0={\mathcal{P}}_\Omega^-[S_{ \lambda },E^\sharp[\Delta u+\lambda u]](x)
}
\\ \nonumber
&& 
+\int_{\partial\Omega}
u(y)\frac{\partial}{\partial\nu_{\Omega,y} }\left(S_{n,\lambda} (x-y)\right)\,d\sigma_{y}
-\langle \frac{\partial u}{\partial\nu_{\Omega}}(y),S_{n,\lambda} (x-y)\rangle 
\quad\forall x\in\Omega^-\,.
\end{eqnarray}
\item[(ii)] If $\Delta u+\lambda u=0$ in $\Omega$, then
\begin{eqnarray}
\label{thm:dthirdgreenh2}
\lefteqn{
\int_{\partial\Omega}
u(y)\frac{\partial}{\partial\nu_{\Omega,y} }S_{n,\lambda} (x-y)
\,d\sigma_{y}-\langle \frac{\partial u}{\partial\nu_{\Omega}}(y),S_{n,\lambda} (x-y)\rangle 
}
\\ \nonumber
&&\qquad\qquad\qquad\qquad\qquad\qquad
=
\left\{
\begin{array}{ll}
u(x) & {\text{if}}\ x\in\Omega\,,
\\
0 & {\text{if}}\ x\in {\mathbb{R}}^{n}\setminus\overline{\Omega}\,.
\end{array}\right.
\end{eqnarray}
\end{enumerate}
\end{theorem}
{\bf Proof.} (i) By the approximation  Lemma \ref{lem:apr0ad}, there exists a sequence $\{u_j\}_{j\in {\mathbb{N}}}$  in $C^{1,\alpha}(\overline{\Omega})_\Delta$ as in (\ref{lem:apr1ad1}). Since $u_j\in C^{1,\alpha}(\overline{\Omega})_\Delta
\subseteq C^2(\Omega)\cap C^{1,\alpha}(\overline{\Omega})$ and $\Delta u_j\in C^{0,\alpha}(\overline{\Omega})$, then the classical third Green Identity (cf.~\textit{e.g.}, \cite[Thm.~9.3 of the Appendix]{La25a}),  equalities  (\ref{prop:nschext3}), (\ref{lem:conoderdeducl1})  and 
\[
{\mathcal{P}}_\Omega[S_{n, \lambda },E^\sharp[\Delta u+\lambda u_j]](x)
=\int_{\Omega}(\Delta u_j(y)+\lambda u_j(y))S_{n,\lambda} (x-y)\,dy
\qquad\forall x\in {\mathbb{R}}^n
\]
imply that
\begin{eqnarray}
\label{thm:dthirdgreenh3}
\lefteqn{
u_j(x)={\mathcal{P}}_\Omega^+[S_{n,\lambda },E^\sharp[\Delta u_j+\lambda u_j]](x)
}
\\ \nonumber
&&  
+\int_{\partial\Omega}
u_j(y)\frac{\partial}{\partial\nu_{\Omega,y} }\left(S_{n,\lambda} (x-y)\right)\,d\sigma_{y}
-\langle \frac{\partial u_j}{\partial\nu_{\Omega}}(y),S_{n,\lambda} (x-y)\rangle 
\quad\forall x\in\Omega\,,
\\ \label{thm:dthirdgreenh3a}
\lefteqn{
0={\mathcal{P}}_\Omega^-[S_{n, \lambda },E^\sharp[\Delta u_j+\lambda u_j]](x)
}
\\ \nonumber
&& 
+\int_{\partial\Omega}
u_j(y)\frac{\partial}{\partial\nu_{\Omega,y} }\left(S_{n,\lambda} (x-y)\right)\,d\sigma_{y}
-\langle \frac{\partial u_j}{\partial\nu_{\Omega}}(y),S_{n,\lambda} (x-y)\rangle 
\quad\forall x\in\Omega^-\,.
\end{eqnarray}
Now let $r\in]0,+\infty[$ be such that $\overline{\Omega}\subseteq 
{\mathbb{B}}(0,r)$ and $\beta\in]0,\alpha[$. By Proposition \ref{prop:dvpsnecr-1a},   we have
\begin{eqnarray}\label{thm:dthirdgreenh4}
&&{\mathcal{P}}_\Omega^+[S_{n,\lambda},E^\sharp[\cdot]]\in 
{\mathcal{L}}\left(C^{-1,\beta}(\overline{\Omega}),C^{1,\beta}(\overline{\Omega})\right)\,,
\\	\nonumber
&&{\mathcal{P}}_\Omega^-[S_{n,\lambda},E^\sharp[\cdot]]_{|\overline{{\mathbb{B}}(0,r)}\setminus\Omega}\in 
{\mathcal{L}}\left(C^{-1,\beta}(\overline{\Omega}),C^{1,\beta}(\overline{{\mathbb{B}}(0,r)}\setminus\Omega)\right)\,.
\end{eqnarray}
Then classical properties of the double layer potential imply that
\begin{eqnarray}\label{thm:dthirdgreenh5}
&&w_\Omega^+[S_{n,\lambda},\cdot]\in 
{\mathcal{L}}\left(C^{0,\alpha}(\partial\Omega),C^{0,\alpha}(\overline{\Omega})\right)\,,
\\	\nonumber
&&w_\Omega^-[S_{n,\lambda},\cdot]_{|\overline{{\mathbb{B}}(0,r)}\setminus\Omega}\in 
{\mathcal{L}}\left(C^{0,\alpha}(\partial\Omega),C^{0,\alpha}(\overline{{\mathbb{B}}(0,r)}\setminus\Omega)\right)\,,
\end{eqnarray}
(cf.~\textit{e.g.}, Theorem \ref{thm:dlay}). Since $S_{n,\lambda}(x-\cdot)\in C^{1,\beta}(\partial\Omega)$, the continuity of the distributional normal derivative from $C^{0,\beta}(\partial\Omega)_\Delta$ to $V^{-1,\beta}(\partial\Omega)$ implies that
\begin{equation}\label{thm:dthirdgreenh6}
\lim_{j\to\infty}\langle \frac{\partial u_j}{\partial\nu_{\Omega}}(y),S_{n,\lambda} (x-y)\rangle 
=\langle \frac{\partial u}{\partial\nu_{\Omega}}(y),S_{n,\lambda} (x-y)\rangle \qquad\forall x\in\Omega\cup\Omega^-\,,
\end{equation}
(cf.~Proposition \ref{prop:ricodnu}). Then by taking the limit as $j$ tends to infinity in equality  (\ref{thm:dthirdgreenh3}) with $x\in\Omega$  and in equality  (\ref{thm:dthirdgreenh3a})     with $x\in {\mathbb{B}}_n(0,r)\setminus\overline{\Omega}$, we deduce the validity of statement (i) for all $x\in\Omega$ and  
$x\in {\mathbb{B}}_n(0,r)\setminus\overline{\Omega}$. Since $r$ such that $\overline{\Omega}\subseteq 
{\mathbb{B}}(0,r)$ has been chosen arbitrarily, we deduce the validity of 
 equality  (\ref{thm:dthirdgreenh1a}) for all $x\in\Omega^-$. Statement (ii) is an immediate consequence of statement (i).\hfill  $\Box$ 

\vspace{\baselineskip}

 Next, we  choose $r\in]0,+\infty[$ such that $\overline{\Omega}\subseteq {\mathbb{B}}_n(0,r)$ and we observe that the linear operator from $C^{1,\alpha}(\partial{\mathbb{B}}_n(0,r))$ to $C^{1,\alpha}((\partial\Omega)\cup(\partial{\mathbb{B}}_n(0,r)))$ that is defined by the equality
\[
\stackrel{o}{E}_{r,\partial\Omega}[v]\equiv\left\{
\begin{array}{ll}
 0 & \text{if}\ x\in \partial\Omega\,,
 \\
 v(x) & \text{if}\ x\in \partial{\mathbb{B}}_n(0,r)\,,
\end{array}
\right.\qquad\forall v\in C^{1,\alpha}(\partial{\mathbb{B}}_n(0,r)) 
\]
is linear and continuous and that accordingly the transpose map $\stackrel{o}{E}_{r,\partial\Omega}^t$ is linear and continuous
\[
\text{from}\ 
\left(C^{1,\alpha}((\partial\Omega)\cup(\partial{\mathbb{B}}_n(0,r)))\right)'\ 
\text{to}\ \left(C^{1,\alpha}(\partial{\mathbb{B}}_n(0,r))\right)'\,.
\]
Then we prove the following immediate consequence of the second Green Identity in the form of Theorem \ref{thm:dsecondgreenh}.
\begin{lemma}\label{lem:orsecgrh}
 Let  $\alpha\in]0,1[$. Let $\Omega$ be a bounded open subset of ${\mathbb{R}}^n$ of class $C^{1,\alpha}$.  Let $\lambda\in {\mathbb{C}}$.  Let $r\in]0,+\infty[$ such that $\overline{\Omega}\subseteq {\mathbb{B}}_n(0,r)$. Let  $u\in C^{0,\alpha}(\overline{{\mathbb{B}}_n(0,r)}\setminus\Omega)$, $v\in C^{1,\alpha}(\overline{{\mathbb{B}}_n(0,r)}\setminus\Omega)_\Delta$. If
\begin{equation}\label{lem:orsecgrh1}
 \Delta u +\lambda u=0
 \quad\text{in}\  {\mathbb{B}}_n(0,r)\setminus\overline{\Omega}\,,
\end{equation}
then
\begin{eqnarray}\label{lem:orsecgrh2}
\lefteqn{
\langle 
\stackrel{o}{E}_{r,\partial\Omega}^t\left[
\partial_{\nu_{{\mathbb{B}}_n(0,r)\setminus\overline{\Omega}}}u
\right],v_{|\partial {\mathbb{B}}_n(0,r)} \rangle-\int_{\partial {\mathbb{B}}_n(0,r)}
\frac{\partial v }{\partial\nu_{{\mathbb{B}}_n(0,r)}}
u\,d\sigma
}
\\ \nonumber
&& 
=
-
\langle\frac{\partial u}{\partial\nu_{\Omega^-}},v\rangle
+
\int_{\partial \Omega}
\frac{\partial v }{\partial\nu_{\Omega^-}}
u\,d\sigma
- \lambda \int_{
{\mathbb{B}}_n(0,r)\setminus\overline{\Omega}
}u  v\,dx
-\int_{
{\mathbb{B}}_n(0,r)\setminus\overline{\Omega}
}u\Delta v\,dx\,.
\end{eqnarray} 
\end{lemma}
{\bf Proof.} Since $u$ satisfies equation (\ref{lem:orsecgrh1}), then $\Delta u\in C^{0,\alpha}(\overline{{\mathbb{B}}_n(0,r)}\setminus\Omega)$ and accordingly 
$u\in C^{0,\alpha}(\overline{{\mathbb{B}}_n(0,r)}\setminus\Omega)_\Delta$.
By the second Green Identity in the form of Theorem \ref{thm:dsecondgreenh} in the set ${\mathbb{B}}_n(0,r)\setminus\overline{\Omega}$, we have
\begin{eqnarray*}
\lefteqn{
\langle\frac{\partial u}{\partial\nu_{{\mathbb{B}}_n(0,r)\setminus\overline{\Omega}}},v\rangle
=\int_{
(\partial\Omega)\cup(\partial {\mathbb{B}}_n(0,r))
}
\frac{\partial v}{\partial\nu_{
{\mathbb{B}}_n(0,r)\setminus\overline{\Omega}
}}u
\,d\sigma
}
\\ \nonumber
&&\qquad\quad
+\langle
E^\sharp_{{\mathbb{B}}_n(0,r)\setminus\overline{\Omega}}[\Delta u],v
\rangle-\int_{
{\mathbb{B}}_n(0,r)\setminus\overline{\Omega}
} u\Delta v\,dx\,.
\end{eqnarray*}
Since $\Delta u=-\lambda u$,  equality (\ref{prop:nschext3}) implies that, 
\begin{eqnarray*}
\lefteqn{
\langle\frac{\partial u}{\partial\nu_{{\mathbb{B}}_n(0,r)\setminus\overline{\Omega}}},
 \stackrel{o}{E}_{\partial\Omega,r}[v_{|\partial\Omega}]  
\rangle
}
\\ \nonumber
&&\qquad\quad
+\langle\frac{\partial u}{\partial\nu_{{\mathbb{B}}_n(0,r)\setminus\overline{\Omega}}},
 v-\stackrel{o}{E}_{\partial\Omega,r}[v_{|\partial\Omega}]  
\rangle
\\ \nonumber
&&\qquad
=\int_{
\partial\Omega 
}
\frac{\partial v}{\partial\nu_{
{\mathbb{B}}_n(0,r)\setminus\overline{\Omega}
}}u
\,d\sigma
+\int_{
 \partial {\mathbb{B}}_n(0,r) 
}
\frac{\partial v}{\partial\nu_{
{\mathbb{B}}_n(0,r)\setminus\overline{\Omega}
}}u
\,d\sigma
\\ \nonumber
&&\qquad\quad
-\lambda\int_{
{\mathbb{B}}_n(0,r)\setminus\overline{\Omega}
} u v\,dx-\int_{
{\mathbb{B}}_n(0,r)\setminus\overline{\Omega}
} u\Delta v\,dx
\end{eqnarray*}
(see (\ref{eq:zerom}) for the definition of $\stackrel{o}{E}_{\partial\Omega,r}$). 
Next we note that 
\[
 v-\stackrel{o}{E}_{\partial\Omega,r}[v_{|\partial\Omega}] 
= \stackrel{o}{E}_{r,\partial\Omega }[v_{|\partial {\mathbb{B}}_n(0,r)}] \,.
\]
Then the Definitions \ref{defn:conoderdedu}, \ref{defn:endedr} of normal derivative imply that 
\begin{eqnarray*}
\lefteqn{
\langle\frac{\partial u}{\partial\nu_{\Omega^-}},v_{|\partial\Omega}\rangle
+
\langle 
\stackrel{o}{E}_{r,\partial\Omega}^t\left[
\partial_{\nu_{{\mathbb{B}}_n(0,r)\setminus\overline{\Omega}}}u
\right],v_{|\partial {\mathbb{B}}_n(0,r)} \rangle 
}
\\ \nonumber
&& 
=\int_{
\partial\Omega 
}
\frac{\partial v}{\partial\nu_{
 \Omega^-
}}u
\,d\sigma
+\int_{
 \partial {\mathbb{B}}_n(0,r) 
}
\frac{\partial v}{\partial\nu_{
{\mathbb{B}}_n(0,r) 
}}u
\,d\sigma
\\ \nonumber
&& \quad
-\lambda\int_{
{\mathbb{B}}_n(0,r)\setminus\overline{\Omega}
} u v\,dx-\int_{
{\mathbb{B}}_n(0,r)\setminus\overline{\Omega}
} u\Delta v\,dx\,.
\end{eqnarray*}
  
\hfill  $\Box$ 

\vspace{\baselineskip}

\begin{remark}\label{rem:orsecgrh} Let  $\alpha\in]0,1[$. Let $\Omega$ be a bounded open subset of ${\mathbb{R}}^n$ of class $C^{1,\alpha}$.  Let $\lambda\in {\mathbb{C}}$.   If   $u\in C^{0,\alpha}_{{\mathrm{loc}}	}(\overline{\Omega^-})$ satisfies equation (\ref{lem:orsecgrh1}) in $\Omega^-$, then 
$u\in C^{0,\alpha}(\overline{{\mathbb{B}}_n(0,r)}\setminus\Omega)_\Delta$ for all  $r\in]0,+\infty[$ such that $\overline{\Omega}\subseteq {\mathbb{B}}_n(0,r)$ and  Lemma \ref{lem:c0ademc1a} implies that
$u\in C^{1,\alpha}_{{\mathrm{loc}}	}( \Omega^-)$. Then 
\cite[\S 2
]{La25} implies that  
   the distributional normal derivative of $u$ coincides with the classical normal derivative of $u$ on $\partial {\mathbb{B}}_n(0,r)$   and  formula (\ref{lem:orsecgrh2}) takes the following form
 \begin{eqnarray}\label{rem:orsecgrh3}
\lefteqn{
\int_{\partial {\mathbb{B}}_n(0,r)}\frac{\partial u}{\partial\nu_{{\mathbb{B}}_n(0,r)}} v\,d\sigma
-
\frac{\partial v }{\partial\nu_{{\mathbb{B}}_n(0,r)}}
u\,d\sigma
}
\\ \nonumber
&& 
=
-
\langle\frac{\partial u}{\partial\nu_{\Omega^-}},v\rangle
+
\int_{\partial \Omega}
\frac{\partial v }{\partial\nu_{\Omega^-}}
u\,d\sigma
-\lambda \int_{
{\mathbb{B}}_n(0,r)\setminus\overline{\Omega}
}u  v\,dx
-\int_{
{\mathbb{B}}_n(0,r)\setminus\overline{\Omega}
}u\Delta v\,dx 
\end{eqnarray} 
for all  $v\in C^{1,\alpha}(\overline{{\mathbb{B}}_n(0,r)}\setminus\Omega)_\Delta$ and  $r\in]0,+\infty[$ such that $\overline{\Omega}\subseteq {\mathbb{B}}_n(0,r)$. 
\end{remark}
We are now ready to prove the second Green Identity for solutions of the Helmholtz equation that satisfy the outgoing radiation condition.
\begin{theorem}\label{thm:desecondgreenh}
 Let  $\alpha\in]0,1[$. Let $\Omega$ be a bounded open subset of ${\mathbb{R}}^n$ of class $C^{1,\alpha}$. Let $k\in {\mathbb{C}}\setminus ]-\infty,0]$, ${\mathrm{Im}}\,k\geq 0$. If $u\in C^{0,\alpha}_{{\mathrm{loc}}	}(\overline{\Omega^-})$ and $v\in C^{1,\alpha}_{{\mathrm{loc}}	}(\overline{\Omega^-})_\Delta$ satisfy the outgoing $(k)$-radiation condition and the equations
\begin{eqnarray}
\label{thm:desecondgreenh1}
\Delta u+k^{2}u=0\qquad{\text{in}}\ \Omega^{-}\,,
\\
\nonumber
\Delta v+k^{2}v=0\qquad{\text{in}}\ \Omega^{-}\,,
\end{eqnarray}
then 
\begin{equation}
\label{thm:desecondgreenh2}
\langle \frac{\partial u}{\partial\nu_{\Omega^-}},v\rangle-
 \int_{\partial\Omega}\frac{\partial v}{\partial\nu_{\Omega^-}}u\,d\sigma=0\,.
\end{equation}
\end{theorem}
{\bf Proof.} By Remark \ref{rem:orsecgrh}, $u\in C^{1,\alpha}_{{\mathrm{loc}}	}( \Omega^-)$ and it suffices to show that the integral in the left hand side of 
(\ref{rem:orsecgrh3}) tends to $0$ as $r$ tends to $+\infty$. To do so we note that 
\begin{eqnarray}
\label{thm:desecondgreenh3}
\lefteqn{
\int_{\partial{\mathbb{B}}_{n}(0,r)}  \frac{\partial u}{\partial\nu_{{\mathbb{B}}_{n}(0,r)}}v-  \frac{\partial v}{\partial\nu_{{\mathbb{B}}_{n}(0,r)}}u\,d\sigma
}
\\ \nonumber
&&
=
\int_{\partial{\mathbb{B}}_{n}(0,r)}v \left(\frac{\partial u}{\partial\nu_{{\mathbb{B}}_{n}(0,r)}}-iku\right)\,d\sigma
-
\int_{\partial{\mathbb{B}}_{n}(0,r)}  \left(\frac{\partial v}{\partial\nu_{{\mathbb{B}}_{n}(0,r)}}-ikv\right)u\,d\sigma\,,
\end{eqnarray}
 that 
\begin{eqnarray}
\label{thm:desecondgreenh4}
\lefteqn{
\left|\int_{\partial{\mathbb{B}}_{n}(0,r)}u \left(\frac{\partial v}{\partial\nu_{{\mathbb{B}}_{n}(0,r)}}-ikv\right)\,d\sigma\right|
}
\\
\nonumber
&&
\qquad\leq
\left(
\int_{\partial{\mathbb{B}}_{n}(0,r)}|u|^{2}\,d\sigma
\right)^{1/2}
\left(
\int_{\partial{\mathbb{B}}_{n}(0,r)}\left|\frac{\partial v}{\partial\nu_{{\mathbb{B}}_{n}(0,r)}}-ikv\right|^{2}\,d\sigma
\right)^{1/2}
\\
\nonumber
&&
\qquad\leq\left(
\int_{\partial{\mathbb{B}}_{n}(0,r)}|u|^{2}\,d\sigma
\right)^{1/2}
s_{n}^{1/2}\sup_{|x|=r}|x|^{\frac{n-1}{2}}\left|\frac{\partial v}{\partial\nu_{{\mathbb{B}}_{n}(0,r)}}-ikv\right|\,,
\end{eqnarray}
for all 
$r\in]0,+\infty[$ such that $\overline{\Omega}\subseteq {\mathbb{B}}_n(0,r)$
and that the same inequality holds by switching $u$ with $v$. Since both $u$ and $v$ satisfy the Helmholtz equation outside of a ball with radius $r_1\in]0,+\infty[$ such that $\overline{\Omega}\subseteq {\mathbb{B}}_n(0,r_1)$ and are of class $C^2$ in ${\mathbb{R}}^n\setminus {\mathbb{B}}_n(0,r_1)$, it is classically known that
\begin{equation}\label{thm:desecondgreenh5}
\limsup_{r\to\infty}\int_{\partial{\mathbb{B}}_{n}(0,r)}
|u|^{2}\,d\sigma<+\infty\,,\quad
\limsup_{r\to\infty}\int_{\partial{\mathbb{B}}_{n}(0,r)}
|v|^{2}\,d\sigma<+\infty\,,
\end{equation}
(cf.~\textit{e.g.}, Colton and Kress \cite[(3.8) in the proof of Theorem 3.3]{CoKr92}). Then by equality (\ref{thm:desecondgreenh3}), by taking the limit as $r$ tends to infinity in the right hand side of (\ref{thm:desecondgreenh4}) and in the inequality that is obtained by switching $u$ with $v$ and by exploiting the outgoing $(k)$-radiation condition, we conclude that the integral in the left hand side of 
(\ref{rem:orsecgrh3}) tends to $0$ as $r$ tends to $+\infty$ and thus the proof is complete.\hfill  $\Box$ 

\vspace{\baselineskip}

Then we have the following distributional  form of the third Green Identity for solutions of the Helmholtz equation in exterior domains.
\begin{theorem}\label{thm:dethirdgreenh}
Let  $\alpha\in]0,1[$. Let $\Omega$ be a bounded open subset of ${\mathbb{R}}^n$ of class $C^{1,\alpha}$. Let $k\in {\mathbb{C}}\setminus ]-\infty,0]$, ${\mathrm{Im}}\,k\geq 0$.  If $u\in C^{0,\alpha}_{ {\mathrm{loc}} }(\overline{\Omega^-}) $ satisfies the outgoing $(k)$-radiation condition and   $\Delta u+k^2 u=0$ in $\Omega^-$, then
\begin{eqnarray}
\label{thm:dethirdgreenh1}
\lefteqn{
\int_{\partial\Omega}
u(y)\frac{\partial}{\partial\nu_{\Omega^-,y} }\left(\tilde{S}_{n,k;r}(x-y)\right)
\,d\sigma_{y}-\langle \frac{\partial u}{\partial\nu_{\Omega^-}}(y),\tilde{S}_{n,k;r} (x-y)\rangle 
}
\\ \nonumber
&&\qquad\qquad\qquad\qquad\qquad\qquad\qquad\qquad\qquad
=
\left\{
\begin{array}{ll}
0 & {\text{if}}\ x\in  \Omega
\,,
\\
 u(x) & {\text{if}}\ x\in\Omega^-\,,
\end{array}\right.
\end{eqnarray}
where $\frac{\partial u}{\partial\nu_{\Omega^-}}$ denotes the distributional normal derivative of Definition \ref{defn:endedr}. 
 \end{theorem}
{\bf Proof.} Let $x\in {\mathbb{R}}^n\setminus\partial\Omega$. Let $r\in]0,+\infty[$ be such that
$\overline{\Omega}\cup\{x\}\subseteq {\mathbb{B}}_n(0,r)$. By the third Green Identity in the form of Theorem \ref{thm:dthirdgreenh} in the set ${\mathbb{B}}_n(0,r)\setminus\overline{\Omega}$, we have
\begin{eqnarray}\label{thm:dethirdgreenh2}
\lefteqn{
 \int_{\partial\Omega\cup\partial {\mathbb{B}}_n(0,r)}u(y)\frac{\partial}{\partial\nu_{
 {\mathbb{B}}_n(0,r)\setminus\overline{\Omega}
 ,y} }\left(\tilde{S}_{n,k;r}(x-y)\right)\,d\sigma
}
\\ \nonumber
&&\qquad\qquad\qquad\qquad
-\langle\partial_{\nu_{{\mathbb{B}}_n(0,r)\setminus\overline{\Omega}}}u(y),\tilde{S}_{n,k;r}(x-y)\rangle=\left\{
\begin{array}{ll}
0 & \text{if}\ x\in\Omega \,,
\\
u(x) & \text{if}\  x\in\Omega^-  \,.
\end{array}
\right.
\end{eqnarray}
Next we note that 
\[
 v-\stackrel{o}{E}_{\partial\Omega,r}[v_{|\partial\Omega}]  = \stackrel{o}{E}_{r,\partial\Omega}[v_{|\partial {\mathbb{B}}_n(0,r)}]   \,.
\]
Then the Definition \ref{defn:endedr}  of normal derivative implies that 
\begin{eqnarray}\label{thm:dethirdgreenh3}
\lefteqn{
 \int_{\partial\Omega}u(y)\frac{\partial}{\partial\nu_{
 \Omega^-
 ,y} }\left(\tilde{S}_{n,k;r}(x-y)\right)\,d\sigma
 }
\\ \nonumber
&&\qquad\quad
+
  \int_{\partial {\mathbb{B}}_n(0,r)}u(y)\frac{\partial}{\partial\nu_{
 {\mathbb{B}}_n(0,r) 
 ,y} }\left(\tilde{S}_{n,k;r}(x-y)\right)\,d\sigma
\\ \nonumber
&&\qquad\quad
-\langle \stackrel{o}{E}_{\partial\Omega,r}^t\partial_{\nu_{ {\mathbb{B}}_n(0,r)\setminus\overline{\Omega}}}u(y),\tilde{S}_{n,k;r}(x-y)\rangle
\\ \nonumber
&&\qquad\quad
-\langle \stackrel{o}{E}_{r,\partial\Omega}^t\partial_{\nu_{ {\mathbb{B}}_n(0,r)\setminus\overline{\Omega}}}u(y),\tilde{S}_{n,k;r}(x-y)\rangle
\\ \nonumber
&&\qquad
=
\int_{\partial\Omega}u(y)\frac{\partial}{\partial\nu_{
 \Omega^-
 ,y} }\left(\tilde{S}_{n,k;r}(x-y)\right)\,d\sigma
 -\langle \partial_{\nu_{\Omega^-} }u(y),\tilde{S}_{n,k;r}(x-y)\rangle
\\ \nonumber
&&\qquad\quad
+\int_{
 \partial {\mathbb{B}}_n(0,r) 
}u(y)\frac{\partial}{\partial\nu_{
 {\mathbb{B}}_n(0,r) 
 ,y} }\left(\tilde{S}_{n,k;r}(x-y)\right)
 \\ \nonumber
&&\qquad\quad
-\langle \stackrel{o}{E}_{r,\partial\Omega}^t\partial_{\nu_{ {\mathbb{B}}_n(0,r)\setminus\overline{\Omega}}}u(y),\tilde{S}_{n,k;r}(x-y)\rangle
 \\ \nonumber
&&\qquad
=\left\{
\begin{array}{ll}
0 & \text{if}\  x\in\Omega
 \,,
\\
  u(x) & \text{if}\ x\in\Omega^-  \,.
\end{array}
\right.
\end{eqnarray}
Since $u\in C^2(\Omega^-)$, 
\cite[\S 2
]{La25} 
implies that 
 the distributional normal derivative of $u$ on $\partial{\mathbb{B}}_n(0,r)$ coincides with the classical normal derivative of $u$ on $\partial {\mathbb{B}}_n(0,r)$ and thus we have
\begin{eqnarray}\label{thm:dethirdgreenh4}
\lefteqn{
\int_{
 \partial {\mathbb{B}}_n(0,r)} 
 u(y)\frac{\partial}{\partial\nu_{
 {\mathbb{B}}_n(0,r) \setminus\overline{\Omega}
 ,y} }\left(\tilde{S}_{n,k;r}(x-y)\right)
 }
  \\ \nonumber
&&\qquad\quad
-\langle \stackrel{o}{E}_{r,\partial\Omega}^t\partial_{\nu_{ 
{\mathbb{B}}_n(0,r)\setminus\overline{\Omega}
}}u(y),\tilde{S}_{n,k;r}(x-y)\rangle
\\ \nonumber
&&\qquad
=
-\int_{
 \partial {\mathbb{B}}_n(0,r)
}u(y)\frac{\partial}{\partial\nu_{
 {\mathbb{B}}_n(0,r)^- 
 ,y} }\left(\tilde{S}_{n,k;r}(x-y)\right)
 \\ \nonumber
&&\qquad\quad
+\int_{\partial {\mathbb{B}}_n(0,r)} 
\frac{\partial u}{\partial\nu_{
 {\mathbb{B}}_n(0,r)^-}}(y),\tilde{S}_{n,k;r}(x-y)\,d\sigma\,.
\end{eqnarray}
Since $u$ and $\tilde{S}_{n,k;r}(x-\cdot)$ belong to $C^2(\overline{{\mathbb{B}}_n(0,r)^-})$ and satisfy the outgoing $(k)$-radiation condition, the second Green Identity for classical solutions of the Helm\-holtz equation that satisfy the outgoing $(k)$-radiation condition implies that the right hand-side of (\ref{thm:dethirdgreenh4}) equals $0$ (cf.~\textit{e.g.}, Theorem \ref{thm:desecondgreenh}). Then equality (\ref{thm:dethirdgreenh3}) implies the validity of the equality of the statement.\hfill  $\Box$ 

\vspace{\baselineskip}

\section{A representation theorem for the interior Neumann eigenfunctions}\label{sec:rinefh}

 We first introduce the following elementary remark. 
\begin{remark}
 Let   $m\in \{0,1\}$, $\alpha\in]0,1[$,  $\lambda\in {\mathbb{C}}$. 	 	Let $\Omega$ be a bounded open subset of ${\mathbb{R}}^{n}$ of class $C^{\max\{1,m\},\alpha}$. If $u\in C^{m,\alpha}(\overline{\Omega})$ satisfies equation $\Delta u+ \lambda 	u=0$ in $\Omega$, then $u\in C^{m,\alpha}(\overline{\Omega})_\Delta$.
\end{remark}
Next we show the following representation theorem for the interior Neumann eigenfunctions. Here we note that the normal derivative is to be interpreted in the sense of Definition \ref{defn:conoderdedu}.
\begin{theorem}\label{thm:eipn}
Let  $m\in {\mathbb{N}}$, $\alpha\in]0,1[$. Let $\Omega$ be a bounded open subset of ${\mathbb{R}}^{n}$ of class $C^{\max\{1,m\},\alpha}$. Then the following statements hold.
\begin{enumerate}
\item[(i)] Let  $\lambda\in {\mathbb{C}}$. If $S_{n,\lambda}$ is a fundamental solution of $\Delta+\lambda$ and if $u\in C^{ m,\alpha}(\overline{\Omega})$ satisfies the problem
\begin{equation}
\label{thm:eipn1}
\left\{
\begin{array}{ll}
\Delta u+\lambda u=0 &{\text{in}}\ \Omega\,,
\\
\frac{\partial u}{\partial\nu_{\Omega}}=0
&{\text{on}}\ \partial\Omega\,,
\end{array}
\right.
\end{equation}
then  $u\in C^{\max\{1,m\},\alpha}(\overline{\Omega})$,  
\begin{equation}
\label{thm:eipn2}
u(x)=  \int_{\partial\Omega}\frac{\partial }{\partial\nu_{\Omega,y}}\left(S_{n,\lambda}(x-y)\right)u
(y)\,d\sigma_{y}\qquad\forall x\in \Omega\,,
\end{equation}
and 
\begin{equation}
\label{thm:eipn3}
-\frac{1}{2}u(x)
+
\int_{\partial\Omega}\frac{\partial }{\partial\nu_{\Omega,y}}\left(S_{n,\lambda}(x-y)\right)u(y)\,d\sigma_{y}=0\qquad\forall x\in\partial\Omega\,.
\end{equation}
\item[(ii)] Let $k\in {\mathbb{C}}\setminus ]-\infty,0]$, ${\mathrm{Im}}\,k\geq 0$. If for each $j\in {\mathbb{N}}\setminus\{0\}$ such that $j\leq \varkappa^{-}$,  
$k^{2}$ is not a Dirichlet eigenvalue for $-\Delta$ in $\Omega_{j}^{-}$, and if $\omega\in C^{m,\alpha}(\partial\Omega)$ satisfies the equation
\begin{equation}
\label{thm:eipn4}
-\frac{1}{2}\omega(x)
+
\int_{\partial\Omega}\frac{\partial }{\partial\nu_{\Omega,y}}\left(\tilde{S}_{n,k;r}(x-y)\right)
\omega(y)\,d\sigma_{y}=0\qquad\forall x\in\partial\Omega\,,
\end{equation}
then the function $u$ from $\overline{\Omega}$ to ${\mathbb{C}}$ defined by 
\begin{equation}
\label{thm:eipn5}
u(x)=w_\Omega^{+}[\tilde{S}_{n,k;r},\omega](x)
\qquad\forall x\in\overline{\Omega}\,,
\end{equation}
belongs to $C^{\max\{1,m\},\alpha}(\overline{\Omega})$ and 
satisfies the Neumann problem (\ref{thm:eipn1}) with $\lambda=k^{2}$,  and 
\begin{equation}
\label{thm:eipn6}
u_{|\partial\Omega}=\omega\qquad{\text{on}}\ \partial
\Omega\,.
\end{equation}
In particular,  $\omega\in C^{\max\{1,m\},\alpha}(\partial\Omega) $. 		 
Moreover, 
$w_\Omega^{-}[\tilde{S}_{n,k;r},\omega]=0$ in $\overline{\Omega^{-}}$.  (For the definition of $\varkappa^-$, see right after (\ref{eq:exto})). 
\end{enumerate}
\end{theorem}
{\bf Proof.}  (i) The membership of $u$ in $C^{\max\{1,m\},\alpha}(\overline{\Omega})$ is obvious in case $m\geq 1$ and is a consequence of the regularization theorem of \cite[\S 4]{La25} in case $m=0$.	 		 If $u$ satisfies the interior Neumann problem (\ref{thm:eipn1}), then
equality (\ref{thm:eipn2}) is an immediate consequence  of the third Green identity   of Theorem \ref{thm:dthirdgreenh}. Then by the continuity of $u$ and of $w_\Omega^{+}[S_{n,\lambda},u_{|\partial\Omega}]$ on $\overline{\Omega}$, we also have
\[
u(x)=w_\Omega^{+}[S_{n,\lambda},u](x)\qquad\forall x\in 
\overline{\Omega}\,,
\]
and accordingly
\[
u(x)=\frac{1}{2}u(x)
+
\int_{\partial\Omega}\frac{\partial }{\partial\nu_{\Omega,y}}\left(S_{n,\lambda}(x-y)\right)
u(y)\,d\sigma_{y}\qquad\forall x\in\partial\Omega\,,
\]
and thus equality (\ref{thm:eipn3}) follows (cf.~Theorem \ref{thm:dlay}). 

We now prove statement (ii). We set $u^{\pm}\equiv w_\Omega^{\pm}[\tilde{S}_{n,k;r},\omega]$. Then
it is classically known that $u^+\in C^{m,\alpha}(\overline{\Omega})$ and that $u^-\in C^{m,\alpha}_{ {\mathrm{loc}} }(\overline{\Omega^{-}})$ (cf.~\textit{e.g.}, Theorem \ref{thm:dlay}, \cite[Thm.~7.3]{DoLa17}).  By equality (\ref{thm:eipn4}) and the jump formula (\ref{thm:dlay1a}) for the acoustic double layer potential, we have  $u^{-}=0$ on $\partial\Omega$  and $u^{+}=\omega$ on $\partial\Omega$.  	Then $u^-$ solves the Dirichlet problem for the Helmholtz equation in $\Omega^{-}_{j}$  and our assumption implies that  $u^{-}=0$ on $\Omega^{-}_{j}$ for each $j\in {\mathbb{N}}\setminus\{0\}$ such that $j\leq \varkappa^{-}$. On the other hand, $u^{-}$ satisfies the outgoing $(k)$-radiation condition and thus $u^{-}=0$ on $\Omega^{-}_{0}$ (cf.~Theorem \cite[
\S 5]{La25}). Accordingly, $u^{-}$ vanishes identically in $\overline{\Omega^-} $ and thus $\frac{\partial u^{-}}{\partial\nu_{\Omega^-}}=0$ on $\partial\Omega^{-}$. Then the jump formula for the normal derivative of the acoustic double layer potential implies that 
$\frac{\partial u^{+}}{\partial\nu_{\Omega}}=-\frac{\partial u^{-}}{\partial\nu_{\Omega^-}}=0$ on $\partial\Omega$ (cf. \cite[\S 6]{La25a}), and thus $u^{+}$ satisfies the Neumann problem (\ref{thm:eipn1}).  Then statement (i) implies that $u\in C^{\max\{1,m\},\alpha}(\overline{\Omega})$.	 	\hfill  $\Box$ 

\vspace{\baselineskip}

\begin{corollary}
\label{corol:eipnn}
  Let  $m\in{\mathbb{N}}$, $\alpha\in]0,1[$. Let $\Omega$ be a bounded open subset of ${\mathbb{R}}^{n}$ of class $C^{\max\{1,m\},\alpha}$. Then the following statements hold.  
\begin{enumerate}
\item[(i)]   Let  $\lambda\in {\mathbb{C}}$. Let $S_{n,\lambda}$ is a fundamental solution of $\Delta+\lambda$. The restriction  operator on the boundary $\partial\Omega$ is a linear injection from the space 
\begin{eqnarray}
\label{corol:eipnn1}
\lefteqn{
\{
u\in C^{m,\alpha}(\overline{\Omega}):\,\Delta u+k^{2}u=0, \frac{\partial u}{\partial\nu_{\Omega}}=0\}
}
\\ \nonumber
&&\qquad
=\{
u\in C^{\max\{1,m\},\alpha}(\overline{\Omega}):\,\Delta u+k^{2}u=0, \frac{\partial u}{\partial\nu_{\Omega}}=0\}
\end{eqnarray}
into the space
\begin{equation}
\label{corol:eipnn2}
\{
\omega\in C^{m,\alpha}(\partial\Omega):\,
-\frac{1}{2}\omega 
+W_\Omega[S_{n,\lambda},\omega]
=0
\}\,.
\end{equation}
Moreover,
\begin{eqnarray}\label{corol:eipnn3}
\lefteqn{
\{
\omega\in C^{m,\alpha}(\partial\Omega):\,
-\frac{1}{2}\omega 
+W_\Omega[S_{n,\lambda},\omega]
=0
\}
}
\\ \nonumber
&&\qquad\qquad 
=\{
\omega\in C^{\max\{1,m\},\alpha}(\partial\Omega):\,
-\frac{1}{2}\omega 
+W_\Omega[S_{n,\lambda},\omega]
=0
\}\,.
\end{eqnarray}
\item[(ii)] Let $k\in {\mathbb{C}}\setminus ]-\infty,0]$, ${\mathrm{Im}}\,k\geq 0$. Assume that for each $j\in {\mathbb{N}}\setminus\{0\}$ such that $j\leq \varkappa^{-}$,  
$k^{2}$ is not a Dirichlet eigenvalue for $-\Delta$ in $\Omega_{j}^{-}$.  Then the restriction  operator on the boundary $\partial\Omega$ is an isomorphism from the space in (\ref{corol:eipnn1})
onto the space
\begin{equation}
\label{corol:eipnn4}
\{
\omega\in C^{\max\{1,m\},\alpha}(\partial\Omega):\,
-\frac{1}{2}\omega 
+W_\Omega[\tilde{S}_{n,k;r},\omega]
=0
\}
\,.
\end{equation}
The  inverse of such an isomorphism is the map $w^{+}_\Omega[\tilde{S}_{n,k;r},\cdot] $. In particular, the dimension of the space in (\ref{corol:eipnn4}) equals the geometric multiplicity of $k^{2}$ as a Neumann eigenvalue of $-\Delta$ in $\Omega$. 
\end{enumerate}
\end{corollary}
{\bf Proof.} (i) By Theorem \ref{thm:eipn}, the  restriction is a linear operator from the space in (\ref{corol:eipnn1}) to the space in (\ref{corol:eipnn2}). By the Holmgren Uniqueness Theorem, the restriction is injective in the space in (\ref{corol:eipnn1}). In order to show the equality in (\ref{corol:eipnn3}), it suffices to assume that $m=0$, that $\omega\in C^{0,\alpha}(\partial\Omega)$ solves
equation 
\begin{equation}\label{corol:eipnn5}
-\frac{1}{2}\omega 
+W_\Omega[S_{n,\lambda},\omega]
=0
\end{equation}
 and to show that $\omega\in C^{1,\alpha}(\partial\Omega)$. If $\omega\in C^{0,\alpha}(\partial\Omega)$ solves
equation (\ref{corol:eipnn5}), then the jump formula (\ref{thm:dlay1a}) for the acoustic double layer potential implies that
\begin{equation}\label{corol:eipnn6}
w_\Omega^-[S_{n,\lambda},\omega]_{|\partial\Omega}=0\,.
\end{equation}
By Theorem \ref{thm:dlay}, we have $w_\Omega^-[S_{n,\lambda},\omega]\in C^{0,\alpha}_{{\mathrm{loc}}}(\overline{\Omega^-})$. Then  equality (\ref{corol:eipnn6}) and the regularization Theorem of \cite[\S 4]{La25} for solutions of the Helmholtz equation implies that $
w_\Omega^-[\tilde{S}_{n,k;r},\omega]\in C^{1,\alpha}_{{\mathrm{loc}}}(\overline{\Omega^-})
$.  Then the jump formula for the (distributional) normal derivative of the double layer potential of \cite[\S 6]{La25a} implies that 
\[
\frac{\partial w_\Omega^+[\tilde{S}_{n,k;r},\omega]}{\partial\nu_{\Omega}}=-\frac{\partial w_\Omega^-[\tilde{S}_{n,k;r},\omega]}{\partial\nu_{\Omega^-}}\in C^{0,\alpha}(\partial\Omega) 
\]
 and the regularization Theorem of \cite[\S 4]{La25} for solutions of the Helmholtz equation implies that $w_\Omega^+[\tilde{S}_{n,k;r},\omega]\in C^{1,\alpha}(\overline{\Omega})$. Then the jump formula for the double layer potential implies that
\[
\omega=w_\Omega^+[\tilde{S}_{n,k;r},\omega]_{|\partial\Omega}-w_\Omega^-[\tilde{S}_{n,k;r},\omega]_{|\partial\Omega}\in C^{1,\alpha}(\partial\Omega)\,.
\]
 Statement (ii) is an immediate consequence of statement (i) and of Theorem \ref{thm:eipn} (ii).
 \hfill  $\Box$ 

\vspace{\baselineskip}

\section{A nonvariational form of the interior Neumann problem for the Helmholtz equation}\label{sec:inpddh}

We now consider the interior (nonvariational) Neumann problem in the case in which the Neumann datum is in the space $V^{-1,\alpha}(\partial\Omega)$. For the corresponding classical result in case $\varkappa^-=0$, we refer to  Colton and Kress \cite[Thm.~3.20]{CoKr92}.  (For the definition of $\varkappa^-$, see right after (\ref{eq:exto})).  
\begin{theorem}[of existence for the interior Neumann problem]\label{thm:exintneuhe}
  Let \,  $\alpha\in]0,1[$. Let $\Omega$ be a bounded open subset of ${\mathbb{R}}^{n}$ of class $C^{1,\alpha}$. 
  
  Let $k\in {\mathbb{C}}\setminus ]-\infty,0]$, ${\mathrm{Im}}\,k\geq 0$. Assume that for each $j\in {\mathbb{N}}\setminus\{0\}$ such that $j\leq \varkappa^{-}$,  
$k^{2}$ is not a Dirichlet eigenvalue for $-\Delta$ in $\Omega_{j}^{-}$.  If  $g\in V^{-1,\alpha}(\partial\Omega)$, then the interior Neumann problem
  \begin{equation}
\label{thm:exintneuhe1}
\left\{
\begin{array}{ll}
\Delta u+k^2 u=0 &{\text{in}}\ \Omega\,,
\\
\frac{\partial u}{\partial\nu_\Omega}=g
&{\text{on}}\ \partial\Omega\,
 \end{array}
\right.
\end{equation}
has   a  solution $u\in C^{0,\alpha}(\overline{\Omega})_\Delta$ if and only if 
\begin{equation}
\label{thm:exintneuhe2}
\langle g,v_{|\partial\Omega}\rangle =0
\end{equation}
for all $v\in C^{1,\alpha}(\overline{\Omega})$ such that
\begin{equation}\label{thm:exintneuhe3}
\Delta v+k^2 v=0 \qquad\text{in}\ \Omega\,,\qquad \frac{\partial v}{\partial\nu_\Omega}=0\qquad\text{on}\ \partial\Omega\,.
\end{equation}
If condition (\ref{thm:exintneuhe2}) is satisfied, then the integral equation
\begin{equation}\label{thm:exintneuhe4}
-\frac{1}{2}\phi+W_\Omega^t[\tilde{S}_{n,k;r},\phi]=g\qquad\text{on}\ \partial\Omega
\end{equation}
has a solution $\phi\in V^{-1,\alpha}(\partial\Omega)$ and the function $u\equiv v_\Omega^+[\tilde{S}_{n,k;r},\phi]$ belongs to $ C^{0,\alpha}(\overline{\Omega})_\Delta$ and solves the interior Neumann problem 
(\ref{thm:exintneuhe2}). 
\end{theorem}
{\bf Proof.}  Assume that a solution $u\in C^{0,\alpha}(\overline{\Omega})_\Delta$ exists and that $v\in C^{1,\alpha}(\overline{\Omega})$ satisfies (\ref{thm:exintneuhe3}). Then the second Green Identity in distributional form of Theorem \ref{thm:dsecondgreenh} and equality (\ref{prop:nschext3}) imply that
\begin{eqnarray*}
\lefteqn{
\langle g,v\rangle =\langle  \frac{\partial u}{\partial\nu_\Omega}, v_{|\partial\Omega}\rangle 
}
\\ \nonumber
&&\qquad\qquad
=\int_{\partial\Omega}u \frac{\partial v}{\partial\nu_\Omega}  \,d\sigma
-
\int_\Omega u(\Delta v+k^2 v)\,dx+\langle E^\sharp[\Delta u+k^2 u],v\rangle =0\,.
\end{eqnarray*}
Conversely, we now assume that the compatibility condition (\ref{thm:exintneuhe2})  holds true. By  Corollary  \ref{corol:eipnn} with $m=1$, we have
\[
\langle g,\psi\rangle =0\qquad\forall\psi\in 
{\mathrm{Ker}}\left(
-\frac{1}{2}I+W_\Omega[\tilde{S}_{n,k;r},\cdot]_{|C^{1,\alpha}(\partial\Omega)}
\right)\,.
\]
By classical results, $W_\Omega[\tilde{S}_{n,k;r},\cdot]$ is compact in $C^{1,\alpha}(\partial\Omega)$ (cf.~\textit{e.g.}, \cite[Cor.~9.1]{DoLa17}) and by  \cite[\S 8]{La25a}, $W_\Omega^t[\tilde{S}_{n,k;r},\cdot]$ is compact in $V^{-1,\alpha}(\partial\Omega)$. Then  the Fredholm Alternative Theorem of Wendland \cite{We67}, \cite{We70} in the duality pairing
\[
\left(V^{-1,\alpha}(\partial\Omega),C^{1,\alpha} (\partial\Omega)\right) 
\]
 implies that there exists $\phi\in  V^{-1,\alpha}(\partial\Omega)$ such that
\[
-\frac{1}{2}\phi+W_\Omega^t[\tilde{S}_{n,k;r},\phi]=g\qquad\text{on}\ \partial\Omega
\]
(cf.~\textit{e.g.}, \cite[Thm.~5.8]{DaLaMu21}).
Thus Theorem \ref{thm:slhco-1a} ensures that $u\equiv v_\Omega^+[\tilde{S}_{n,k;r},\phi]$ belongs to $C^{0,\alpha}(\overline{\Omega})_\Delta$ and the jump formulas of \cite[\S 6]{La25a} ensures that $u$ satisfies the boundary condition of the interior Neumann problem (\ref{thm:exintneuhe1}).
Hence, $u$ solves the interior Neumann problem (\ref{thm:exintneuhe1}).\hfill  $\Box$ 

\vspace{\baselineskip}

\section{A representation theorem for the exterior Neumann eigenfunctions}\label{sec:renefh}
We first introduce the following elementary remark. 
\begin{remark}
 Let   $m\in \{0,1\}$, $\alpha\in]0,1[$, $\lambda\in {\mathbb{C}}$. 	  Let $\Omega$ be a bounded open subset of ${\mathbb{R}}^{n}$ of class $C^{\max\{1,m\},\alpha}$. If  $u\in C^{m,\alpha}_{{\mathrm{loc}}}(\overline{\Omega^-})$ satisfies equation $\Delta u+ \lambda	 	u=0$ in $\Omega^-$, then $u\in C^{m,\alpha}_{{\mathrm{loc}}}(\overline{\Omega^-})_\Delta$.
\end{remark}
 Next we show the following representation theorem for the exterior Neumann eigenfunctions. Here we note that the normal derivative is to be interpreted in the sense of Definition \ref{defn:endedr}.
\begin{theorem}
\label{thm:eipen}
Let   $m\in {\mathbb{N}}$, $\alpha\in]0,1[$. Let $\Omega$ be a bounded open subset of ${\mathbb{R}}^{n}$ of class $C^{\max\{1,m\},\alpha}$. Let $k\in {\mathbb{C}}\setminus ]-\infty,0]$, ${\mathrm{Im}}\,k\geq 0$. Then the following statements hold.
\begin{enumerate}
\item[(i)] If $u\in C^{ m,\alpha}_{{\mathrm{loc}} }(\overline{\Omega^-})  $ satisfies the problem
\begin{equation}
\label{thm:eipen1}
\left\{
\begin{array}{ll}
\Delta u+k^2 u=0 &{\text{in}}\ \Omega^-\,,
\\
-\frac{\partial u}{\partial\nu_{\Omega^-}}=0
&{\text{on}}\ \partial\Omega\,,
\\
u\ \text{satisfies\ the\ outgoing}\  (k)-\text{radiation\ condition}\,,
\end{array}
\right.
\end{equation}
then  $u\in C^{\max\{1,m\},\alpha}_{{\mathrm{loc}}}(\overline{\Omega^-})$,  
\begin{equation}
\label{thm:eipen2}
u(x)=  -\int_{\partial\Omega}\frac{\partial }{\partial\nu_{\Omega,y}}\left(\tilde{S}_{n,k;r}(x-y)\right)u
(y)\,d\sigma_{y}\qquad\forall x\in \Omega^-\,,
\end{equation}
and 
\begin{equation}
\label{thm:eipen3}
\frac{1}{2}u(x)
+
\int_{\partial\Omega}\frac{\partial }{\partial\nu_{\Omega,y}}\left(\tilde{S}_{n,k;r}(x-y)\right)u(y)\,d\sigma_{y}=0\qquad\forall x\in\partial\Omega\,.
\end{equation}
\item[(ii)]   If  $k^{2}$ is not a Dirichlet eigenvalue for $-\Delta$ in $\Omega$, and if $\omega\in C^{m,\alpha}(\partial\Omega)$ satisfies the equation
\begin{equation}
\label{thm:eipen4}
\frac{1}{2}\omega(x)
+
\int_{\partial\Omega}\frac{\partial }{\partial\nu_{\Omega,y}}\left(\tilde{S}_{n,k;r}(x-y)\right)
\omega(y)\,d\sigma_{y}=0\qquad\forall x\in\partial\Omega\,,
\end{equation}
then the function $u$ from $\overline{\Omega^-}$ to ${\mathbb{C}}$ defined by 
\begin{equation}
\label{thm:eipen5}
u(x)=w_\Omega^{-}[\tilde{S}_{n,k;r},-\omega](x)
\qquad\forall x\in\overline{\Omega^-}\,,
\end{equation}
belongs to $C^{ \max\{1,m\},\alpha}_{{\mathrm{loc}} }(\overline{\Omega^-}) $ and 
satisfies the Neumann problem (\ref{thm:eipen1})  and
\begin{equation}
\label{thm:eipen6}
w_\Omega^{-}[\tilde{S}_{n,k;r},-\omega]_{|\partial\Omega}=\omega\qquad{\text{on}}\ \partial
\Omega\,.
\end{equation}
In particular,  $\omega\in C^{\max\{1,m\},\alpha}(\partial\Omega) $. 		 
Moreover, 
$w_\Omega^{+}[\tilde{S}_{n,k;r},-\omega]=0$ in $\overline{\Omega}$. 
\end{enumerate}
\end{theorem}
{\bf Proof.}  (i) The membership of $u$ in $C^{ \max\{1,m\},\alpha}_{{\mathrm{loc}} }(\overline{\Omega^-}) $ is obvious in case $m\geq 1$ and is a consequence of the regularization theorem of \cite[\S 4]{La25} in case $m=0$.	 	
If $u$ satisfies the exterior Neumann problem (\ref{thm:eipen1}), then the third Green Identity for exterior domains of Theorem \ref{thm:dethirdgreenh} implies the validity of equality (\ref{thm:eipen2}). Then by the continuity of $u$ and of $w^{-}_\Omega[\tilde{S}_{n,k;r},u]$ on $\overline{\Omega^-}$, we also have
\[
u(x)=-w^{-}_\Omega[\tilde{S}_{n,k;r},u](x)\qquad\forall x\in 
\overline{\Omega^-}\,,
\]
(cf.~Theorem \ref{thm:dlay}). Then the jump formula (\ref{thm:dlay1a}) for the double layer acoustic potential implies that
\[
u(x)=\frac{1}{2}u(x)
-
\int_{\partial\Omega}\frac{\partial }{\partial\nu_{\Omega}(y)}\left(\tilde{S}_{n,k;r}(x-y)\right)
u(y)\,d\sigma_{y}\qquad\forall x\in\partial\Omega\,,
\]
and thus equality (\ref{thm:eipen3}) follows. We now prove statement (ii). We set $u^{\pm}\equiv -w^{\pm}_\Omega[\tilde{S}_{n,k;r},\omega]$. Then
it is classically known that $u^+\in C^{m,\alpha}(\overline{\Omega})$ and that $u^-\in C^{m,\alpha}_{ {\mathrm{loc}} }(\overline{\Omega^{-}})$ (cf.~\textit{e.g.}, Theorem \ref{thm:dlay}, \cite[Thm.~7.3]{DoLa17}). By the jump formula (\ref{thm:dlay1a}) for the double layer acoustic potential  and by equality (\ref{thm:eipen4}), we have  $u^{+}=0$ on $\partial\Omega$ and
\[
-\omega=w^{+}_\Omega[\tilde{S}_{n,k;r},-\omega]-w^{-}_\Omega[\tilde{S}_{n,k;r},-\omega]=-u^-
\quad\text{on}\ \partial\Omega\,.
\]
Then $u^+$ solves the  Dirichlet problem for the Helmholtz equation in $\Omega$ with $0$ boundary datum and our assumption implies that  $u^{+}=0$ in $\overline{\Omega}$. Then $\frac{\partial u^{+}}{\partial\nu_{\Omega}}=0$ on $\partial\Omega$. Then the jump formula for the (distributional) normal derivative of the double layer potential (cf.~\cite[\S 6]{La25a}) implies that 
$-\frac{\partial u^{-}}{\partial\nu_{\Omega^-}}=\frac{\partial u^{+}}{\partial\nu_{\Omega}}=0$ on $\partial\Omega$ and thus $u^{-}$ satisfies the exterior Neumann problem (\ref{thm:eipen1}).  Then statement (i) implies that $u\in C^{ \max\{1,m\}	,\alpha}_{{\mathrm{loc}} }(\overline{\Omega^-}) $.	 	\hfill  $\Box$ 

\vspace{\baselineskip}

\begin{corollary}
\label{corol:eipenn}
Let  $m\in{\mathbb{N}}$, $\alpha\in]0,1[$. Let $\Omega$ be a bounded open subset of ${\mathbb{R}}^{n}$ of class $C^{\max\{1,m\},\alpha}$. Let $k\in {\mathbb{C}}\setminus ]-\infty,0]$, ${\mathrm{Im}}\,k\geq 0$. Then the following statements hold.  
\begin{enumerate}
\item[(i)]    The restriction  operator on the boundary $\partial\Omega$ is a linear injection  from the space 
\begin{eqnarray}
\label{corol:eipenn1}
\lefteqn{\{
u\in C^{m,\alpha}_{{\mathrm{loc}}}(\overline{\Omega^-}):\,\Delta u+k^{2}u=0 \ \text{in}\ \Omega^-, \frac{\partial u}{\partial\nu_{\Omega^-}}=0\ \text{on}\ \partial\Omega\,,
}
\\ \nonumber
&&\quad\qquad
u \ \text{satisfies\ the\ outgoing\ $(k)$-radiation\ condition} \}
\\ \nonumber
&&\quad\ 
=\{
u\in C^{\max\{1,m\},\alpha}_{{\mathrm{loc}}}(\overline{\Omega^-}):\,\Delta u+k^{2}u=0 \ \text{in}\ \Omega^-, \frac{\partial u}{\partial\nu_{\Omega^-}}=0\ \text{on}\ \partial\Omega\,,
\\ \nonumber
&&\quad\qquad
u \ \text{satisfies\ the\ outgoing\ $(k)$-radiation\ condition} \}
\end{eqnarray}
to the space
\begin{equation}
\label{corol:eipenn2}
\{
\omega\in C^{m,\alpha}(\partial\Omega^-):\,\frac{1}{2}\omega 
+W_\Omega[\tilde{S}_{n,k;r},\omega]
=0\}\,.
\end{equation}
Moreover,
\begin{eqnarray}\label{corol:eipenn3}
\lefteqn{\{
\omega\in C^{m,\alpha}(\partial\Omega^-):\,\frac{1}{2}\omega 
+W_\Omega[\tilde{S}_{n,k;r},\omega]
=0\}
 }
\\ \nonumber
&&\qquad\qquad
=\{
\omega\in C^{\max\{1,m\},\alpha}(\partial\Omega^-):\,\frac{1}{2}\omega 
+W_\Omega[\tilde{S}_{n,k;r},\omega]
=0\}
\,.
\end{eqnarray}

\item[(ii)]     Assume  that  $k^{2}$ is not a Dirichlet eigenvalue for $-\Delta$ in $\Omega$. 
 Then the restriction  operator on the boundary $\partial\Omega$ is an isomorphism from the space 
in (\ref{corol:eipenn1}) onto the space
\begin{equation}
\label{corol:eipenn4}
\{
\omega\in C^{\max\{1,m\},\alpha}(\partial\Omega):\,
\frac{1}{2}\omega 
+W_\Omega[\tilde{S}_{n,k;r},\omega]
=0\}\,.
\end{equation}
The  inverse of such an isomorphism is the map $-w^{-}_\Omega[\tilde{S}_{n,k;r},\cdot] $. In particular, the dimension of the space in (\ref{corol:eipenn4}) equals the geometric multiplicity of $k^{2}$ as a Neumann eigenvalue of $-\Delta$ in $\Omega^-$.  
\end{enumerate}
\end{corollary}  
{\bf Proof.}  (i) By Theorem \ref{thm:eipen}, the  restriction is a linear operator from the space in (\ref{corol:eipenn1}) to the space in (\ref{corol:eipenn2}). By the Holmgren Uniqueness Theorem, the restriction is injective in the space in (\ref{corol:eipenn1}). In order to show the equality in (\ref{corol:eipenn3}), it suffices to assume that $m=0$, that $\omega\in C^{0,\alpha}(\partial\Omega)$ solves
equation 
\begin{equation}\label{corol:eipenn5}
\frac{1}{2}\omega 
+W_\Omega[S_{n,\lambda},\omega]
=0
\end{equation}
 and to show that $\omega\in C^{1,\alpha}(\partial\Omega)$. If $\omega\in C^{0,\alpha}(\partial\Omega)$ solves
equation (\ref{corol:eipenn5}), then the jump formula (\ref{thm:dlay1a}) for the acoustic double layer potential implies that
\begin{equation}\label{corol:eipenn6}
w_\Omega[S_{n,\lambda},\omega]_{|\partial\Omega}=0\,.
\end{equation}
By Theorem \ref{thm:dlay}, we have $w_\Omega[S_{n,\lambda},\omega]\in C^{0,\alpha}(\overline{\Omega})$. Then  equality (\ref{corol:eipenn6}) and the regularization Theorem of \cite[\S 4]{La25} for solutions of the Helmholtz equation implies that $w_\Omega[\tilde{S}_{n,k;r},\omega]\in C^{1,\alpha}(\overline{\Omega})$. Then the jump formula for the (distributional) normal derivative of the double layer potential (cf.~\cite[\S 6]{La25a}) implies that 
\[
-\frac{\partial w_\Omega^-[\tilde{S}_{n,k;r},\omega]}{\partial\nu_{\Omega^-}}
=
\frac{\partial w_\Omega^+[\tilde{S}_{n,k;r},\omega]}{\partial\nu_{\Omega}}\in C^{0,\alpha}(\partial\Omega) 
\]
and the regularization Theorem of \cite[\S 4]{La25} for solutions of the Helmholtz equation implies that $w_\Omega^-[\tilde{S}_{n,k;r},\omega]\in C^{1,\alpha}_{{\mathrm{loc}}}(\overline{\Omega^-})$. Then the jump formula for the double layer potential implies that
\[
\omega=w_\Omega^+[\tilde{S}_{n,k;r},\omega]_{|\partial\Omega}-w_\Omega^-[\tilde{S}_{n,k;r},\omega]_{|\partial\Omega}\in C^{1,\alpha}(\partial\Omega)\,.
\]
 Statement (ii) is an immediate consequence of statement (i) and of Theorem \ref{thm:eipen} (ii).

\hfill  $\Box$ 

\vspace{\baselineskip}

\begin{remark}
Under the assumptions of Corollary \ref{corol:eipenn} (ii),  the uniqueness
of Neumann eigenfunctions that satisfy the  outgoing $(k)$-radiation condition in $(\Omega^{-})_{0}$ (cf.~\cite[\S 5]{La25}) implies that the dimension of the space in (\ref{corol:eipenn4})  equals the sum of the geometric multiplicities of $k^{2}$ as a Neumann eigenvalue of $-\Delta$ in $(\Omega^{-})_{j}$ for $j\in {\mathbb{N}}\setminus\{0\}$ and $j\leq \varkappa^{-}$.  (For the definition of $\varkappa^-$, see right after (\ref{eq:exto})).

In particular, if we assume that for  each $j\in {\mathbb{N}}\setminus\{0\}$ such that $j\leq \varkappa^{-}$, $k^2$ is not a Neumann eigenvalue for $-\Delta$ in $(\Omega^{-})_{j}$, then the dimension of the space in (\ref{corol:eipenn2})  equals $0$. For example, this is the case if we assume  that $\Omega^-$ is connected.
 \end{remark}

\section{A nonvariational form of the exterior Neumann problem for the Helmholtz equation}\label{sec:enpddh}

We now consider the exterior (nonvariational) Neumann problem in the case in which the Neumann datum is in the space $V^{-1,\alpha}(\partial\Omega)$. For the corresponding classical result in case $\Omega^-$ is connected and the Neumann data satisfy some extra regularity  assumption,  we refer to  Colton and Kress \cite[Thm.~3.25]{CoKr92}.
\begin{theorem}[of existence for the exterior Neumann problem]\label{thm:exextneuhe}
  Let \,   $\alpha\in]0,1[$. Let $\Omega$ be a bounded open subset of ${\mathbb{R}}^{n}$ of class $C^{1,\alpha}$. 
  
  Let $k\in {\mathbb{C}}\setminus ]-\infty,0]$, ${\mathrm{Im}}\,k\geq 0$. Assume that $k^{2}$ is not a Dirichlet eigenvalue for $-\Delta$ in $\Omega$.  If  $g\in V^{-1,\alpha}(\partial\Omega)$, then the exterior Neumann problem
  \begin{equation}
\label{thm:exextneuhe1}
\left\{
\begin{array}{ll}
\Delta u+k^2 u=0 &{\text{in}}\ \Omega^-\,,
\\
-\frac{\partial u}{\partial\nu_{\Omega^-}}=g
&{\text{on}}\ \partial\Omega\,,
\\
u\ \text{satisfies\ the\ outgoing}\  (k)-\text{radiation\ condition}\,,
 \end{array}
\right.
\end{equation}
has   a  solution $u\in C^{0,\alpha}_{{\mathrm{loc}}}(\overline{\Omega^-})_\Delta$ if and only if 
\begin{equation}
\label{thm:exextneuhe2}
\langle g,v_{|\partial\Omega}\rangle =0
\end{equation}
for all $v\in C^{1,\alpha}_{{\mathrm{loc}}}(\overline{\Omega^-})$ such that
\begin{eqnarray}\label{thm:exextneuhe3}
&&\Delta v+k^2 v=0 \qquad\text{in}\ \Omega^-\,,\qquad \frac{\partial v}{\partial\nu_{\Omega^-}}=0\qquad\text{on}\ \partial\Omega\,,
\\\nonumber
&&v \ \text{satisfies\ the\ outgoing\ $(k)$-radiation\ condition.}
\end{eqnarray}
If condition (\ref{thm:exextneuhe2}) is satisfied, then the integral equation
\begin{equation}\label{thm:exextneuhe4}
\frac{1}{2}\phi+W_\Omega^t[\tilde{S}_{n,k;r},\phi]=g\qquad\text{on}\ \partial\Omega
\end{equation}
has a solution $\phi\in V^{-1,\alpha}(\partial\Omega)$ and the function $u\equiv v_\Omega^-[\tilde{S}_{n,k;r},\phi]$ belongs to $C^{0,\alpha}_{{\mathrm{loc}}}(\overline{\Omega^-})_\Delta$ and solves the exterior Neumann problem 
(\ref{thm:exextneuhe2}). 
\end{theorem}
{\bf Proof.} Assume that  $u\in C^{0,\alpha}_{{\mathrm{loc}}}(\overline{\Omega^-})_\Delta$  solves the exterior Neumann problem 
(\ref{thm:exextneuhe2}) and that $v\in C^{1,\alpha}_{{\mathrm{loc}}}(\overline{\Omega^-})$ satisfies  (\ref{thm:exextneuhe3}). Then the second Green Identity for exterior domains of Theorem \ref{thm:desecondgreenh} implies that
\[
\langle g,v_{|\partial\Omega}\rangle=\langle-\frac{\partial u}{\partial\nu_{\Omega^-}},v_{|\partial\Omega}\rangle
=-\int_{\partial\Omega}u\frac{\partial v}{\partial\nu_{\Omega^-}}\,d\sigma=0\,.
\]
Conversely, we now assume that the compatibility condition (\ref{thm:exextneuhe2})  holds true.  By  Corollary  \ref{corol:eipenn} with $m=1$, we have
\begin{eqnarray*}
\lefteqn{
\biggl\{v_{|\partial\Omega}:\, v\in C^{1,\alpha}_{{\mathrm{loc}}}(\overline{\Omega^-})\,,\ 
\Delta v+k^{2}v=0 \ \text{in}\ \Omega^-, \frac{\partial v}{\partial\nu_{\Omega^-}}=0\ \text{on}\ \partial\Omega\,,
}
\\ \nonumber
&&\qquad\qquad\qquad\qquad\qquad\qquad 
v \ \text{satisfies\ the\ outgoing\ $(k)$-radiation\ condition} \biggr\}
\\ \nonumber
&&\qquad
=\{
\omega\in C^{1,\alpha}(\partial\Omega^-):\,(\ref{thm:eipen4}) \ \text{holds}\}\,.
\end{eqnarray*}
Hence, 
\[
\langle g,\psi\rangle =0\qquad\forall\psi\in 
{\mathrm{Ker}}\left(
\frac{1}{2}I+W_\Omega[\tilde{S}_{n,k;r},\cdot]_{|C^{1,\alpha}(\partial\Omega)}
\right)\,.
\]
By classical results, $W_\Omega[\tilde{S}_{n,k;r},\cdot]$ is compact in $C^{1,\alpha}(\partial\Omega)$ (cf.~\textit{e.g.}, \cite[Cor.~9.1]{DoLa17}) and by  \cite[\S 8]{La25a}, $W_\Omega^t[\tilde{S}_{n,k;r},\cdot]$ is compact in $V^{-1,\alpha}(\partial\Omega)$. Then  the Fredholm Alternative Theorem of Wendland in the duality pairing
\[
\left(V^{-1,\alpha}(\partial\Omega),C^{1,\alpha} (\partial\Omega)\right) 
\]
 implies that there exists $\phi\in  V^{-1,\alpha}(\partial\Omega)$ such that
\[
\frac{1}{2}\phi+W_\Omega^t[\tilde{S}_{n,k;r},\phi]=g\qquad\text{on}\ \partial\Omega\,,
\]
(cf.~\textit{e.g.}, \cite[Thm.~5.8]{DaLaMu21}). Then Theorem \ref{thm:slhco-1a} ensures that $u\equiv v_\Omega^-[\tilde{S}_{n,k;r},\phi]$ belongs to $C^{0,\alpha}_{{\mathrm{loc}}}(\overline{\Omega^-})_\Delta$, Theorem \ref{thm:pora} ensures that $u$ satisfies  the outgoing $(k)$-radiation condition and the jump formulas of \cite[\S 6]{La25a} ensures that $u$ satisfies the boundary condition of the exterior Neumann problem (\ref{thm:exextneuhe2}). Hence, $u$ solves the exterior Neumann problem (\ref{thm:exextneuhe2}).\hfill  $\Box$ 

\vspace{\baselineskip}

  \noindent
{\bf Statements and Declarations}\\

 \noindent
{\bf Competing interests:} This paper does not have any  conflict of interest or competing interest.

 \noindent
{\bf Acknowledgement.}  

The author  acknowledges  the support of   GNAMPA-INdAM  and   of the Project funded by the European Union – Next Generation EU under the National Recovery and Resilience Plan (NRRP), Mission 4 Component 2 Investment 1.1 - Call for tender PRIN 2022 No. 104 of February, 2 2022 of Italian Ministry of University and Research; Project code: 2022SENJZ3 (subject area: PE - Physical Sciences and Engineering) ``Perturbation problems and asymptotics for elliptic differential equations: variational and potential theoretic methods''.

\end{document}